\documentclass{svmult}

\usepackage{amsmath}
\usepackage{amssymb}
\usepackage{latexsym}

\usepackage{euscript}

\usepackage{natbib}

\bibpunct{[}{]}{}{}{}{}

\newtheorem{df}{Definition}[section]
\newtheorem{thm}[df]{Theorem}       
\newtheorem{prop}[df]{Proposition}
\newtheorem{cor}[df]{Corollary}
\newtheorem{ex}[df]{Example}
\newtheorem{exs}[df]{Examples}

\newtheorem{rmk}[df]{Remark}

\newcommand{\pf}{\noindent{\sc Proof.}\ }

\newcommand{\ked}{\hspace*{\fill}%
\lower3pt\hbox{\vrule height1.1ex width .9ex depth -.2ex}%
\vskip9pt}

\newcommand{\quend}{\hspace*{\fill}\lower3pt\hbox{$\boxtimes$}} 

\renewcommand{\mathcal}[1]{\EuScript{#1}}

\renewcommand{\phi}{\varphi}
\newcommand{\chigh}{{\raise1.5pt\hbox{$\chi$}}}

\newcommand{\tilD}{\widetilde{D}}
\newcommand{\tilo}{\widetilde 0}           
\newcommand{\tilq}{\widetilde q}

\newcommand{\VBgpds}{$\mathcal{VB}$--groupoids}

\newcommand{\R}{\mathbb{R}}    
\newcommand{\Z}{\mathbb{Z}}

\newcommand{\gof}{\mathfrak{f}}
\newcommand{\goF}{\mathfrak{F}}

\newcommand{\id}{{\rm id}} 

\newcommand{\xclam}{^{\textstyle !}}

\newcommand{\pback}[1]{\mathbin{\times_{#1}}}

\newcommand{\sol}{\bullet}

\newcommand{\llangle}{\langle\!\langle}
\newcommand{\rrangle}{\rangle\!\rangle}

\newcommand{\lhangle}{\,\rule[-0.1ex]{0.4ex}{0.7\baselineskip}\,\,}
\newcommand{\rhangle}{\,\,\rule[-0.1ex]{0.4ex}{0.7\baselineskip}\,}

\newcommand{\ddt}[1]{\left.\displaystyle\frac{d}{dt}#1\right|_0}

\renewcommand{\Bar}{\overline}

\newcommand{\co}{\colon\thinspace} %
\newcommand{\st}{\ \vert\ }
\newcommand{\isom}{\cong}

\renewcommand{\leq}{\leqslant}

\newcommand{\surj}{-\!\!\!-\!\!\!-\!\!\!\gg}
\newcommand{\inj}{>\!\!\!-\!\!\!-\!\!\!-\!\!\!>}

\newcommand{\lrah}%
{\hbox{$\,-\!\!\!-\!\!\!-\!\!\!-\!\!\!-\!\!\!-\!\!\!-\!\!\!\longrightarrow\,$}}

\newcommand{\mlra}{\hbox{$\,-\!\!\!-\!\!\!\longrightarrow\,$}}

\newcommand{\gpd}{\,\lower1pt\hbox{$\longrightarrow$}\hskip-.24in%
\raise2pt\hbox{$\longrightarrow$}\,}

\newcommand{\dpl}{\mbox{$+\hskip-6pt +\hskip4pt$}}
\newcommand{\DPL}[1]{\mbox{$+\hskip-6pt +_{#1}\hskip4pt$}}

\newcommand{\dminus}%
	{\raise2pt\hbox{\vrule height1pt width 2ex}\hskip3pt}
\newcommand{\dtimes}{\mathbin{\hbox{\huge.}}}

\newcommand{\plusA}{\ \lower 5pt%
	\hbox{${\stackrel{\textstyle +}{\scriptscriptstyle A}}$}\ }
\newcommand{\minusA}{\ \lower 5pt%
	\hbox{${\stackrel{\textstyle -}{\scriptscriptstyle A}}$}\ }
\newcommand{\timesA}{\ \lower 4pt%
	\hbox{${\stackrel{\textstyle .}{\scriptscriptstyle A}}$}\ }

\newcommand{\plusEst}{\ \lower 5pt%
	\hbox{${\stackrel{\textstyle +}{\scriptscriptstyle E^*}}$}\ }
\newcommand{\minusEst}{\ \lower 5pt%
	\hbox{${\stackrel{\textstyle -}{\scriptscriptstyle E^*}}$}\ }
\newcommand{\timesEst}{\ \lower 4pt%
	\hbox{${\stackrel{\textstyle .}{\scriptscriptstyle E^*}}$}\ }

\newcommand{\plusB}{\ \lower 5pt%
	\hbox{${\stackrel{\textstyle +}{\scriptscriptstyle B}}$}\ }
\newcommand{\minusB}{\ \lower 5pt%
	\hbox{${\stackrel{\textstyle -}{\scriptscriptstyle B}}$}\ }
\newcommand{\timesB}{\ \lower 4pt%
	\hbox{${\stackrel{\textstyle .}{\scriptscriptstyle B}}$}\ }

\newcommand{\plusCst}{\ \lower 5pt%
	\hbox{${\stackrel{\textstyle +}{\scriptscriptstyle C^*}}$}\ }
\newcommand{\minusCst}{\ \lower 5pt%
	\hbox{${\stackrel{\textstyle -}{\scriptscriptstyle C^*}}$}\ }
\newcommand{\timesCst}{\ \lower 4pt%
	\hbox{${\stackrel{\textstyle .}{\scriptscriptstyle C^*}}$}\ }

\newcommand{\plusW}{\ \lower 5pt%
	\hbox{${\stackrel{\textstyle +}{\scriptscriptstyle W}}$}\ }
\newcommand{\minusW}{\ \lower 5pt%
	\hbox{${\stackrel{\textstyle -}{\scriptscriptstyle W}}$}\ }
\newcommand{\timesW}{\ \lower 4pt%
	\hbox{${\stackrel{\textstyle .}{\scriptscriptstyle W}}$}\ }

\newcommand{\tvb}{\mbox{\rm I\!I\!I}}

\newcommand{\duer}%
{\mbox{${}^\times\kern-.39em\rule{0.3pt}{8pt}\,\,$}}

\newcommand{\du}%
{\mbox{$\times\kern-.41em\rule[-0.5pt]{0.3pt}{6.5pt}\,$}}

\begin{document}

\abovedisplayskip=6pt plus3pt minus3pt
\belowdisplayskip=6pt plus3pt minus3pt

\title*{{\bf DUALITY AND TRIPLE STRUCTURES}}

\author{Kirill C. H. Mackenzie}

\institute{Department of Pure Mathematics\\
University of Sheffield\\
Sheffield S3 7RH\\
United Kingdom\\
{\sf K.Mackenzie@sheffield.ac.uk}}

\maketitle

\begin{flushright}\it
Dedicated to Alan Weinstein\\
on the occasion of his sixtieth birthday
\end{flushright}

\bigskip

\begin{abstract}
We recall the basic theory of double vector bundles and the
canonical pairing of their duals introduced by the author and 
by Konieczna and Urba\'nski. We then show that the relationship 
between a double vector bundle and its two duals can be understood
simply in terms of an associated cotangent triple vector bundle 
structure. In particular we show that the dihedral group of the
triangle acts on this triple via forms of the isomorphisms $R$
introduced by the author and Ping Xu. We then consider the three 
duals of a general triple vector bundle and show that the 
corresponding group is neither the dihedral group of the square, 
nor the symmetry group on four symbols. 
\end{abstract}

Double structures first appeared 
in Poisson geometry with Alan's ground--break\-ing work on 
symplectic groupoids \citep*{CosteDW}, \citep{Weinstein:1987} and
Poisson groupoids \citep{Weinstein:1988}. The most fundamental 
example of a symplectic groupoid, the cotangent groupoid 
$T^*G$ of an arbitrary Lie groupoid $G$, introduced in 
\citep{CosteDW}, is a groupoid object in the category of
vector bundles. An arbitrary Poisson 
Lie group can be integrated to a symplectic double groupoid
\citep{LuW:1989}. At a simpler level, a Poisson structure on 
a vector bundle is linear \citep{Courant:1990} if and only if the 
associated anchor is a morphism of certain double vector bundles. 

All these phenomena involve doubles in the categorical sense:
taking $\mathcal{S}$ to denote, for example, `vector bundle' or 
`Lie groupoid', a double $\mathcal{S}$ is an $\mathcal{S}$
object in the category of all $\mathcal{S}$. (Groupoid objects
in the category of vector bundles, named \emph{\VBgpds} by 
\citet{Pradines:1988}, may be regarded as double groupoids 
of a special type.) More generally, multiple $\mathcal{S}$ 
structures are the $n$--fold extension of this notion of double. 

The key link between Poisson geometry and double structures
lies in properties of the Poisson anchor. If a Poisson manifold
$P$ is a vector bundle on base $M$, then the Poisson structure
is linear if and only if $\pi^\#\co T^*P\to TP$ is a morphism
of double vector bundles. If $P$ is instead a Lie groupoid on
base $M$, then the groupoid is a Poisson groupoid if and only if
$\pi^\#$ is a morphism of \VBgpds. Thus the Poisson anchor
naturally appears as a map of double structures, and indeed
many of the surprising basic features of Poisson and symplectic 
groupoids are not really so much consequences of Poisson or 
symplectic geometry, as consequences of the duality properties of 
the associated double structures. This point of view is developed 
further in \citep{Mackenzie:GT2004}; in particular the theory of 
Poisson groupoids may be developed entirely in terms of groupoid 
theory and double structures of various kinds.  

The present paper is concerned with the duality of double and 
higher multiple vector bundles. A double vector bundle has two 
duals which are themselves in duality and we show here that 
the various combinations of the two dualization operations gives 
rise to the dihedral (or symmetric) group of order six. We show 
in \S\ref{sect:tvbs} and \S\ref{sect:corners} that a double vector 
bundle and its two duals form the three lower faces of a triple 
vector bundle, the opposite vertex of which is the cotangent of 
the double space. This encapsulates and makes symmetric the 
relations between a double vector bundle and its duals, which 
can otherwise seem rather involved. One may think of three double 
vector bundles with a common vertex and appropriate pairings as 
constituting a two (sic) dimensional version of the familiar 
notion of pairing of vector bundles; we call this a cornering. 

In \S\ref{sect:dtvb} we consider the process of dualizing the 
structures in a triple vector bundle. This may appear to be a
routine extension of the double case, but we show that the
group of dualization operations here is not the dihedral group
of the square, or the symmetric group on four symbols, but a
group of order 72. This appears to demonstrate that the behaviour 
of duality for $n$--fold vector bundles may be a less routine
extension of the double case than one might have expected. 
In the final \S\ref{sect:gp}, we formulate some general 
principles which we believe do hold for general multiple vector 
bundles. 

The study of general double vector bundles was begun by
\citet{Pradines:DVB}, though the case of the tangent
double of an ordinary vector bundle (\ref{ex:TE9.1}) had been 
used in connection theory since the late 1950s. More than a decade 
later, \citet{Pradines:1988} introduced the dualization
process for \VBgpds; in the case of double vector bundles this is 
the duality construction presented here in \S\ref{sect:ddvb}. 
Theorem \ref{thm:dualduality} is from \citep{Mackenzie:1999} and 
was also found by \citet{KoniecznaU:1999}. The idea of deriving 
the pairing (\ref{eq:3duals}) from the cotangent triple was noted 
in \citep{Mackenzie:2002}. 

The results of \S\ref{sect:proldual} first appeared in 
the paper \citep{MackenzieX:1994} of Ping Xu and myself. They
are here obtained as a consequence of the general duality of double
vector bundles. An expanded account of the double case may be 
found in \citep{Mackenzie:GT2004}.

I am grateful to Moty Katzman for introducing me to the {\sf GAP}
software, which was very valuable in \S\ref{sect:dtvb}. 

\section{Double vector bundles}
\label{sect:dvbs}

\begin{df}
A \emph{double vector bundle}
$(D;A,B;M)$ is a system of four vector bundle structures
\begin{equation}                             
\label{diag:dvb}
\begin{matrix}
        &&q^D_B&&\\
        &D&\lrah&B&\\
        &&&&\\
        q^D_A&\Big\downarrow&&\Big\downarrow&q_B\\
        &&&&\\
        &A&\lrah&M&\\
        &&q_A&&\\
\end{matrix}
\end{equation}
in which $D$ has two vector bundle structures, on bases $A$ and $B$,
which are themselves vector bundles on $M$,
such that each of the structure maps of each vector bundle
structure on $D$ (the bundle projection, addition, scalar
multiplication and the zero section) is a morphism of vector bundles
with respect to the other structure.
\end{df}

We refer to $A$ and $B$ as the \emph{side bundles} of $D$, and to 
$M$ as the \emph{double base}. In the two side bundles the 
addition, scalar multiplication and subtraction are denoted by the 
usual symbols $+$, juxtaposition, and $-$. We distinguish the two 
zero--sections, writing
$0^A \co M\to A,\ m\mapsto 0^A_m$, and $0^B \co M\to B,\ m\mapsto
0^B_m$. We may denote an element $d\in D$ by $(d;a,b;m)$ to indicate
that $a = q^D_A(d),\ b = q^D_B(d),\ m = q_B(q^D_B(d)) =
q_A(q^D_A(d)).$

The notation $q^D_A$ is clear; when the base of the bundle is the
double base we write $q_A$, for example, rather than $q^A_M$. 

In the \emph{vertical bundle structure} on 
$D$ with base $A$ the vector bundle operations are denoted 
$\plusA, \timesA, \minusA$, with 
$\tilo^A\co A \to D,\ a\mapsto\tilo^A_a$, for the zero--section. 
Similarly, in the \emph{horizontal bundle 
structure} on $D$ with base $B$ we write 
$\plusB, \timesB, \minusB$ and $\tilo^B\co B\to D,\ 
b\mapsto \tilo^B_b$. For $m\in M$ the 
\emph{double zero}
$\rule[-8pt]{0pt}{22pt}\tilo^A_{0^A_m} = \tilo^B_{0^B_m}$ is denoted 
$\odot_m$ or $0^2_m$. The two
structures on $D$, namely $(D,q^D_B, B)$ and $(D,q^D_A, A)$, will
occasionally be denoted $\tilD_B$ and $\tilD_A$, respectively.

In dealing with general double vector bundles such as 
(\ref{diag:dvb}), we thus usually label objects and operations in 
the two structures on $D$ by the symbol for the base over which they 
take place. The words `horizontal' and `vertical' may be used as an 
alternative, but need to be referred to the arrangement in the 
diagram (\ref{diag:dvb}) or the sequence in $(D;A,B;M)$. When 
considering examples in which $A = B$,  the words `horizontal' 
and `vertical' become necessary, and we use $H$ and $V$ as labels to 
distinguish the two structures on $D$. 

Although the concept of double vector bundle is symmetric, most 
examples are not; in the sequel it will be important to distinguish 
between (\ref{diag:dvb}) and its \emph{flip} in Figure
\ref{fig:flip_TE}(a), 
\begin{figure}
\begin{center}
$\begin{matrix}
        &&  q^D_A  &&\\
        & D^{\mbox{{\tiny flip}}}&\longrightarrow&A&\\
        &&&&\\
  q^D_B &\Big\downarrow&&\Big\downarrow&\\
        &&&&\\
        &B&\longrightarrow&M&\\
        &&&&\\
        &&\mbox{(a)}&&\\
\end{matrix}
\hfil
\begin{matrix}
        && T(q)     &&\\
&TE&\longrightarrow&TM&\\
        &&&&\\
        &\Big\downarrow&&\Big\downarrow&\\
        &&&&\\
        &E&\longrightarrow&M&\\
        &&&&\\
        &&\mbox{(b)}&&\\
\end{matrix}
$
\end{center}
\caption{\ \label{fig:flip_TE}}
\end{figure}
in which the arrangement of the two structures is reversed. 
In such processes it is not the absolute arrangement which is 
significant, but the distinction between whichever arrangement is 
taken at the start, and its flip. 

The condition that each addition in $D$ is a morphism with respect
to the other is:
\begin{equation}                               \label{eq:interchange}
(d_1\plusB d_2)\plusA(d_3\plusB d_4)=
(d_1\plusA d_3)\plusB(d_2\plusA d_4)
\end{equation}
for quadruples $d_1, \ldots, d_4\in D$ such that 
$q^D_B(d_1)= q^D_B(d_2),\ q^D_B(d_3)=q^D_B(d_4)$,
$q^D_A(d_1)=q^D_A(d_3)$, and $q^D_A(d_2)=q^D_A(d_4)$.
Next,
\begin{equation}
t\timesA(d_1\plusB d_2)=t\timesA d_1\plusB t\timesA d_2,
\end{equation}
for $t\in \mathbb{R}$ and $d_1, d_2\in D$ with 
$q^D_B(d_1)=q^D_B(d_2)$; similarly
\begin{equation}
t\timesB(d_1\plusA d_2)=t\timesB d_1\plusA t\timesB d_2,
\end{equation}
for $t\in\R$ and $d_1, d_2\in D$ with $q^D_A(d_1)=q^D_A(d_2)$.
The two scalar multiplications are related by
\begin{equation}
t\timesA(u\timesB d)=u\timesB(t\timesA d),
\end{equation}
where $t,u\in\R$ and $d\in D$.

Lastly, for compatible $a,a'\in A$ and compatible $b,b'\in B$, and 
$t\in\R$,
\begin{equation}
\tilo^A_{a+a'} = \tilo^A_a\plusB\tilo^A_{a'},\qquad
\tilo^A_{ta} = t\timesB\tilo^A_a,
\end{equation}
and
\begin{equation}              \label{eq:zeros}
\tilo^B_{b+b'} = \tilo^B_b\plusA\tilo^B_{b'},\qquad
\tilo^B_{tb} = t\timesA\tilo^B_b.
\end{equation}
Equations (\ref{eq:interchange})--(\ref{eq:zeros}) are known as the 
\emph{interchange laws}.

\begin{df}
A \emph{morphism of double vector bundles}
$$
(\phi;\phi_A,\phi_B;f)\co (D;A,B;M)\to (D';A',B';M')
$$
consists of maps $\phi\co D\to D',\ \phi_A\co A\to A'$,
$\phi_B\co B\to B'$, $f\co M\to M',$ such that each of 
$(\phi, \phi_B)$, $(\phi, \phi_A)$, $(\phi_A,f)$ and $(\phi_B,f)$ 
is a morphism of the relevant vector bundles.

If $M = M'$ and $f\ = \id_M$, we say that $\phi$ is \emph{over $M$}; 
if, further, $A = A'$ and $\phi_A = \id_{A}$, we say that $\phi$ is 
\emph{over $A$} or \emph{preserves $A$}. If $A = A'$ and $B = B'$ 
and both $\phi_A$ and $\phi_B$ are identities, we say that $\phi$ 
\emph{preserves the side bundles}.
\end{df}

\begin{ex}\rm
\label{ex:TE9.1}
For an ordinary vector bundle $(E, q, M)$, applying the tangent
functor to each of the bundle operations yields the 
\emph{tangent prolongation vector bundle} $(TE, T(q), TM)$. 
The zero section is $T(0)\co TM\to TE$. We denote the addition 
by $\dpl$ and the
scalar multiplication and subtraction by $\dtimes$ and $\dminus$. 
Together with the standard structure $(TE, p_E, E)$, we have a double
vector bundle $(TE;E,TM;M)$, shown in Figure~\ref{fig:flip_TE}(b). 
There is no preferred arrangement for the side bundles of $TE$. 
\quend
\end{ex}

\begin{ex}\rm                          \label{ex:trivialdvb}
Let $A, B$ and $C$ be any three vector bundles on the one base $M$, 
and write $D$ for the pullback manifold $A\times_M B\times_M C$ over 
$M$. Then $D$ may be regarded as the direct sum 
$q_A\xclam B \oplus q_A\xclam C$ over $A$, and as the 
direct sum $q_B\xclam A\oplus q_B\xclam C$ over $B$,
and with respect to these two structures, $D$ is a double vector 
bundle with side bundles $A$ and $B$. We call this the 
\emph{trivial double vector bundle over $A$ and $B$ with core 
$C$}\index{double vector bundle!trivial}. %
It is tempting, but incorrect, to denote it by 
$A\oplus B\oplus C$.
\quend
\end{ex}

\begin{ex}\rm
A double vector bundle $(D;A,B;M)$ may be pulled back over both
of its side structures simultaneously. Suppose given vector
bundles $(A', q_{A'}, M')$ and $(B', q_{B'}, M')$ and morphisms 
$\phi\co A'\to A$ and $\psi\co B'\to B$, both over a map 
$f\co M'\to M$. Let $D'$ denote the set of all $(a',d,b')$ such
that $\phi(a') = q^D_A(d)$, $\psi(b') = q^D_B(d)$ and
$q_{A'}(a') = q_{B'}(b')$. Then, with the evident structures,
$(D';A',B';M')$ is a double vector bundle and the projection
$D'\to D$ is a morphism over $\phi$, $\psi$ and $f$. 
\quend
\end{ex}

Further examples follow later in the paper. 

\section{The core and the core sequences}
\label{sect:cores}

Until Example \ref{ex:cores}, consider a fixed double vector bundle
$(D;A,B;M)$. Each of the bundle projections is a morphism with
respect to the other structure and so has a kernel (in the ordinary
sense); denote by $C$ the intersection of the two kernels:
$$
C=\{c \in D\st \exists m\in M \mbox{\ such that}\ q^D_B(c )= 0^B_m,\
q^D_A(c)=0^A_m\}.
$$
This is an embedded submanifold of $D$. We will show that it has a
well--defined vector bundle structure with base $M$, projection $q_C$ 
which is the restriction of $q_B\circ q^D_B = q_A\circ q^D_A$ and 
addition and scalar multiplication which are the restrictions 
of either of the operations on $D$.

Note first that the two additions coincide on $C$ since
$$
c \plusB c' = (c \plusA\odot_m)\plusB(\odot_m\plusA c')=
(c \plusB\odot_m)\plusA(\odot_m\plusB c')=c \plusA c',
$$
for $c,c'\in C$ with $q_C(c) = q_C(c')$, using (\ref{eq:interchange}).
From this it follows that $t\timesB c = t\timesA c$ for integers 
$t$, and consequently for rational $t$, and thence for all real 
$t$ by continuity.

It will often be helpful to distinguish between $c\in C$, regarding 
$C$ as a distinct vector bundle, and the image of $c$ in $D$, which 
we will denote by $\Bar{c}$. Given $c,c'\in C$ with 
$q_C(c ) = q_C(c')$ there is a unique $c + c'\in C$ with
$$
\Bar{c +c'} = \Bar{c }\plusB\Bar{c'} = \Bar{c }\plusA\Bar{c'},
$$
and given $t\in\R$ there is a unique $tc \in C$ such that
$$
\Bar{tc } = t\timesB\Bar{c } = t\timesA\Bar{c }.
$$
It is now easy to prove that $(C,q_C,M)$ is a (smooth) vector bundle,
which we call the 
\emph{core}\index{double vector bundle!core}\index{core!of double vector bundle}
of $(D;A,B;M)$.

\begin{thm}                         
\label{thm:coreseq}
There is an exact sequence
\begin{equation}                             
\label{eq:Vcore}
q_A\xclam C \stackrel{\tau_A}{\inj}
\tilD_A\stackrel{(q^D_B)\xclam}{\surj}q_A\xclam B
\end{equation}
of vector bundles over $A$, and an exact sequence
\begin{equation}                             
\label{eq:Hcore}
q_B\xclam C\stackrel{\tau_B}{\inj}
\tilD_B\stackrel{(q^D_A)\xclam}{\surj}q_B\xclam A
\end{equation}
of vector bundles over $B$, where the injections are
$\tau_A \co(a, c)\mapsto \tilo_a^A \plusB \Bar{c}$ and
$\tau_B \co(b, c)\mapsto \tilo_b^B \plusA \Bar{c}$, respectively, and
$(q^D_B)\xclam$ and $(q^D_A)\xclam$ denote the maps induced by $q^D_B$
and $q^D_A$ into the pullback bundles.
\end{thm}

\pf
Take $a\in A_m,\ c \in C_m$ where $m\in M$. Then both
$\tilo^A_a$ and $\Bar{c}$ project under $q^D_B$ to $0^B_m$. So
$\tilo^A_a\plusB\Bar{c}$ is defined and also projects under $q^D_B$ to
$0^B_m$. That $\tau_A$ is linear over $A$ follows from the 
interchange laws.

Suppose that $d\in D$ has $q^D_B(d) = 0^B_m$ for some $m\in M$. Write 
$a = q^D_A(d)$. Then $d \minusB \tilo^A_a$ is defined and 
$q^D_B(d \minusB\tilo^A_a) = 0^B_m$. On the other hand, 
$\tilq^A(d \minusB\tilo^A_a) = a -a = 0^A_m$. So 
$d \minusB \tilo^A_a \in C_m$. This establishes the exactness 
of (\ref{eq:Vcore}). The proof of (\ref{eq:Hcore}) is similar.
\ked

We refer to (\ref{eq:Vcore}) as the 
\emph{core sequence of $D$ over $A$}, 
and to (\ref{eq:Hcore}) as the \emph{core sequence of $D$ over $B$.}

If $(\phi;\phi_A,\phi_B;f)\co (D;A,B;M)\to (D';A',B';M')$ is a
morphism of double vector bundles, then $\phi\co D\to D'$ maps $C$
into $C'$, the core of $D'$. It is clear that the restriction, 
$\phi_C \co C\to C'$, is a morphism of the vector bundle structures 
on the cores, over $f\co M\to M'$.

\begin{exs}\rm                          
\label{ex:cores}
For $E$ an ordinary vector bundle, consider the tangent double 
vector bundle $(TE;E,TM;M)$. The kernel of $T(q)$ consists of the 
vertical tangent vectors and the kernel of $p_E$ consists of the 
vectors tangent along the zero section; their intersection is
naturally identified with $E$ itself. For clarity we distinguish
$X\in E$ from the core element $\Bar{X}\in TE$. 

The injection map for $TE$ over $E$ is the map $\tau$ which sends 
$(X,Y)\in E_m\times E_m$ to the vector in $E_m$ which has tail 
at $X$ and is parallel to $Y$. In terms of the prolongation 
structure, $\tau(X, Y) = \tilo_X\dpl \Bar{Y}$. 
The injection map over $TM$ is $\upsilon\co (x,Y) 
\mapsto T(0)(x) + \Bar{Y}.$

For $\phi\co E\to E'$ a morphism of vector bundles over 
$f\co M\to M'$, the morphism $T(\phi)$ of the tangent double vector 
bundles induces $\phi\co E\to E'$ on the cores. In the case where 
$\phi$ and $f$ are surjective submersions, the vertical subbundles 
form a double vector subbundle (in an obvious sense) 
$(T^\phi E; E, T^fM;M)$ of $TE$, the core of which is the kernel 
(in the ordinary sense) of $\phi$. 

The trivial double vector bundle $A\times_M B\times_M C$ of 
\ref{ex:trivialdvb} has core $C$.
\quend
\end{exs}

\section{Duals of double vector bundles}            
\label{sect:ddvb}

Throughout this section we consider a double vector bundle as in
(\ref{diag:dvb}), with core bundle $C$. We will show that dualizing 
either structure on $D$ leads again to a double vector bundle; in 
the case of the dual of the structure over $A$ we denote this by
\begin{equation}                 \label{diag:E^*}
\begin{matrix}
        &&q_{C^*}^{\du A}&&\\
        &D\duer A&\lrah&C^*&\\
        &&&&\\
q_A^{\du A}&\Big\downarrow& &\Big\downarrow&q_{C^*}\\
        &&&&\\
        &A&\lrah&M,&\\
        &&q_A&&\\
\end{matrix}
\end{equation}
Here $C^*$ is the ordinary dual of $C$ as a vector bundle over $M$. 
We denote the dual of $D$ as a vector bundle over $A$ by 
$D\duer A$. (We will later modify this notation for cases in 
which $A$ and $B$ are identical.)

The vertical structure in (\ref{diag:E^*}) is the usual dual of the 
bundle structure on $D$ with base $A$, and $q_{C^*}\co C^*\to M$ is 
the usual dual of $q_C\co C\to M.$ The additions and scalar 
multiplications in the side bundles of (\ref{diag:E^*}) will be 
denoted by the usual plain symbols as before. In the two 
structures on $D\duer A$ we write $\plusA, \timesA, \minusA$ 
and $\plusCst, \timesCst, \minusCst$. The zero of $D\duer A$ 
above $a\in A$ is denoted $\tilo_a^{\du A}$.

The unfamiliar projection $q^{\du A}_{C^*}\co (D\duer A)\to C^*$ is 
defined by
\begin{equation}                             \label{eq:unfproj}
\langle q^{\du A}_{C^*}(\Phi), c\rangle =
   \langle\Phi, \tilo^A_a \plusB \Bar c\rangle
\end{equation}
where $c\in C_m,\ \Phi\co (q^D_A)^{-1}(a)\to\R$ and $a\in A_m$. The
addition $\plusCst$ in $D^{\duer A}\to C^*$ is defined by
\begin{equation}                  %
\langle\Phi\plusCst\Phi', d\plusB d'\rangle =
   \langle\Phi,d\rangle + \langle\Phi',d'\rangle
\end{equation}
That this is well--defined depends strongly on the condition
$q^{\du A}_{C^*}(\Phi) = q^{\du A}_{C^*}(\Psi).$ Similarly, define
$$
\langle t\timesCst\Phi, t\timesB d\rangle = t\langle\Phi,d\rangle,
$$
for $t\in\R$ and $d\in D$ with $q^D_A(d) = \tilq^{\du A}_A(\Phi)$.

The zero above $\kappa\in C^*_m$ is denoted $\tilo^{\du A}_\kappa$
and is defined by
\begin{equation}
\langle\,\tilo^{\du A}_\kappa, \tilo^B_b \plusA \Bar c\rangle =
\langle\kappa, c\rangle
\end{equation}
where $b\in B_m,\ c\in C_m.$ The core element $\Bar\psi$ 
corresponding to $\psi\in B_m^*$ is
$$
\langle\Bar\psi, \tilo^B_b \plusA \Bar c\rangle =
\langle\psi,c\rangle.
$$

It is straightforward to verify that (\ref{diag:E^*}) is a double 
vector bundle, and that its core is $B^*$. We call (\ref{diag:E^*}) 
the \emph{vertical dual}\index{double vector bundle!vertical dual} 
or \emph{dual over $A$} of (\ref{diag:dvb}).
\index{duality!of double vector bundle}

As for any double vector bundle, there are exact sequences
\begin{equation}           \label{des:A}
q\xclam_A B^*\stackrel{\sigma_A}{\inj} D\duer A
\stackrel{(q^{\du A}_{C^*})\xclam}{\surj} q_A\xclam C^*,
\end{equation}
of vector bundles over $A$ and
\begin{equation}
q_{C^*}\xclam B^*\stackrel{\sigma_{C^*}}{\inj} D\duer A
\stackrel{(q^{\du A}_A)\xclam }{\surj} q_{C^*}\xclam A,
\end{equation}
of vector bundles over $C^*$. Here the two injections are given by
$$
\sigma_A(a,\psi) = \tilo^{\duer A}_a \plusCst \Bar{\psi},\qquad
\sigma_{C^*}(\kappa,\psi) = \tilo^{\duer A}_\kappa \plusA \Bar{\psi},
$$
where $a\in A, \psi\in B^*, \kappa\in C^*.$ It is easily seen that
$$
\langle\sigma_A(a,\psi), d\rangle = \langle\psi, q^D_B(d)\rangle
$$
for $d\in D$ and so $\sigma_A$ is precisely the dual of 
$(q^D_B)\xclam$. It is clear from the definition of 
$q^{\du A}_{C^*}$ that $(q^{\du A}_{C^*})\xclam = \tau_V^*$. Thus 
(\ref{des:A}) is precisely the dual of the core exact sequence 
(\ref{eq:Vcore}).

For the sequence over $C^*$ we have
$$
\langle\sigma_{C^*}(\kappa,\psi), \tilo^B_b \plusA \Bar{c}\rangle
 = \langle\kappa, c\rangle + \langle\psi, b\rangle
$$
for $\kappa\in  C^*_m,\ \psi\in B_m^*,\ x\in B_m,\ c\in C_m.$

The proof of the following result is straightforward. In 
Figure~\ref{fig:dualmorph} and in similar figures in future, we
omit arrows which are the identity.

\begin{prop}                     
\label{prop:dmor}
Consider a morphism of double vector bundles, as in 
Figure~\emph{\ref{fig:dualmorph}(a)}, which preserves the horizontal
side bundles, 
\begin{figure}[htb]
\setlength{\unitlength}{1cm}
\begin{picture}(10,5.5)
\put(0,4){$\begin{matrix}
                      &&   &\\
                      &D&\longrightarrow &B\\
                      &&&\\
                      &\Big\downarrow& &\Big\downarrow\\
                      &&&\\
                      &A&\longrightarrow &M\\
            \end{matrix}
                      $}
\put(0.5,4.5){\vector(2,-1){1.8}}                     
\put(1.6,4.6){\vector(2,-1){1.8}}                     

\put(0.9,4){$\phi$}
\put(3,4){$\phi_B$}

\put(2.1,2){$\begin{matrix}
                     &&    &\\
                     &D'&\longrightarrow &B'\\
                     &&&\\
                     &\Big\downarrow&&\Big\downarrow\\
                     &&&\\
                     &A&\longrightarrow & M\\
                     &&&\\
		     &&&\\
                     &{\rm (a)}&&\\
               \end{matrix}
                     $}


\put(6,4){$\begin{matrix}
                      &&   &\\
                      &D'\duer A&\longrightarrow &(C')^*\\
                      &&&\\
                      &\Big\downarrow& &\Big\downarrow\\
                      &&&\\
                      &A&\longrightarrow &M\\
            \end{matrix}
                      $}
\put(7,4.5){\vector(2,-1){2}}               
\put(8.5,4.6){\vector(2,-1){2}}               

\put(7,3.8){$\phi\duer A$}
\put(9.5,4.3){$\phi_C^*$}

\put(8.8,1.9){$\begin{matrix}
                     &&    &\\
                     &D\duer A&\longrightarrow &C^*\\
                     &&&\\
                     &\Big\downarrow&&\Big\downarrow\\
                     &&&\\
                     &A&\longrightarrow & M\\
                     &&&\\
		     &&&\\
                     &{\rm (b)}&&\\               
                     \end{matrix}
                     $}
\end{picture}
\caption{\ \label{fig:dualmorph}}
\end{figure}
and which has core morphism $\phi_C\co C\to C'$, where $C'$ is 
the core of $D'$. Dualizing $\phi$ as a morphism of vector bundles 
over $A$, we obtain a morphism $\phi\duer A$ of double vector 
bundles 
over $A$ and $\phi_C^*$, as in Figure~\emph{\ref{fig:dualmorph}(b)}, 
with core morphism $\phi_B^*$. 
\end{prop}

This completes the description of the vertical dual of 
(\ref{diag:dvb}). There is of course also a \emph{horizontal dual}
\begin{equation}                    \label{eq:hdual}
\begin{matrix}
        &&q^{\du B}_B&&\\
        &D\duer B&\lrah&B&\\
        &&&&\\
q^{\du B}_{C^*}&\Big\downarrow&&\Big\downarrow&q_B\\
        &&&&\\
        &C^*&\lrah&M,&\\
        &&q_{C^*}&&\\
\end{matrix}
\end{equation}
with core $A^*\to M$, defined in an analogous way. 

The following result is an entirely new phenomenon, arising from the
double structures. 

\begin{thm}{\rm \citep{Mackenzie:1999}, \citep{KoniecznaU:1999}}
\label{thm:dualduality}
There is a natural (up to sign) duality between the bundles 
$D\duer A$ and $D\duer B$ over $C^*$ given by
\begin{equation}                       
\label{eq:3duals}
\lhangle\Phi, \Psi\rhangle = 
\langle\Phi, d\rangle_A - \langle\Psi, d\rangle_B
\end{equation}
where $\Phi\in D\duer A,\ \Psi\in D\duer B$ have
$q^{\du A}_{C^*}(\Phi) = q^{\du B}_{C^*}(\Psi)$ and $d$ 
is any element of $D$ with $q^D_A(d) = q^{\du A}_A(\Phi)$ and
$q^D_B(d) = q^{\du B}_B(\Psi).$
\end{thm}

Each of the pairings on the RHS of (\ref{eq:3duals}) is a canonical 
pairing of an ordinary vector bundle with its dual, the subscripts 
there indicating the base over which the pairing takes place. 

\bigskip

\pf
Let $\Phi$ and $\Psi$ have the forms $(\Phi; a, \kappa; m)$ and
$(\Psi; \kappa, b; m)$. Then $d$ must have the form $(d; a, b; m)$.
If $d'$ also has the form $(d'; a, b; m)$ then there is a $c\in C_m$
such that $d = d' \plusA (\tilo^A_a \plusB \Bar c)$, and so
$$
\langle\Phi,d\rangle_A = 
\langle\Phi, d'\rangle_A + \langle\kappa, c\rangle
$$
by (\ref{eq:unfproj}). By the interchange law (\ref{eq:interchange})
we also have $d = d' \plusB (\tilo^B_b \plusA \Bar c)$ and so
$$
\langle\Psi,d\rangle_B = \langle\Psi, d'\rangle_B + \langle\kappa, c\rangle.
$$
Thus (\ref{eq:3duals}) is well defined. To check that it is bilinear
is routine. It remains to prove that it is non--degenerate. 

Suppose $\Phi$, given as above, is such that 
$\lhangle\Phi, \Psi\rhangle = 0$ for all 
$\Psi\in(q^{\du B}_{C^*})^{-1}(\kappa).$ Take any $\phi\in A_m^*$
and consider $\Psi = \tilo_\kappa^{\du B}\plusB\Bar\phi$. Then, taking
$d = \tilo^A_a$ we find $\langle\Phi, d\rangle_A = 0$ and
$\langle\Psi, d\rangle_B = \langle\phi, a\rangle.$
Thus $\langle\phi, a\rangle = 0$ for all $\phi\in A_m^*$ and so
$a = 0^A_m$. It therefore follows from the horizontal exact 
sequence for $D\duer A$ that
$$
\Phi = \tilo_\kappa^{\du A} \plusA\Bar\psi
$$
for some $\psi\in B_m^*$. Now taking any $c\in C_m$ and defining
$d = \tilo_b^B\plusA \Bar c$, we find that
$$
\langle\Phi, d\rangle_A = \langle\kappa, c\rangle + \langle\psi, b\rangle
\qquad\mbox{and}\qquad
\langle\Psi, d\rangle_B = \langle\kappa, c\rangle.
$$
So $\langle\psi, b \rangle = 0$ for all $b\in B_m$, since a suitable
$\Psi$ exists for any given $b$. It follows that $\psi = 0\in B_m^*$
and so $\Phi$ is indeed the zero element over $\kappa$. Thus the pairing
(\ref{eq:3duals}) is nondegenerate.
\ked

Note several special cases:
\begin{gather}
\lhangle\tilo^{\du A}_\kappa,\tilo^{\du B}_\kappa\rhangle = 0,\qquad
\lhangle\tilo^{\du A}_a,\tilo^{\du B}_b\rhangle = 0,\\
\lhangle\tilo^{\du A}_a, \Bar\phi\rhangle 
	= -\langle\phi, a\rangle,\qquad
\lhangle\Bar\psi, \tilo_b^{\du B}\rhangle 
	= \langle\psi, b\rangle,\qquad
\lhangle\Bar\psi, \Bar\phi\rhangle = 0, \\
\lhangle\tilo_\kappa^{\du A}\plusA\Bar{\psi}, \Psi\rhangle 
	= \langle\psi, b\rangle,\qquad\qquad
\lhangle\Phi, \tilo^{\du B}_\kappa\plusB\Bar\phi\rhangle 
	= -\langle\phi, a\rangle,
\end{gather}
where $b\in B,\ a\in A,\ \phi\in A^*,\ \psi\in B^*$ and we have
$(\Psi;\kappa, b; m)\in D\duer B$ and 
$(\Phi;a, \kappa; m)\in D\duer A$.

Although we have proved that $D\duer A$ and $D\duer B$ are dual as 
vector bundles over $C^*$, we have not yet considered the 
relationships between the other structures present. This is taken 
care of by the following result, the proof of which is 
straightforward.

\begin{prop}{\rm \citep{Mackenzie:1999}}
\label{prop:dudvb}
Let $(D;A,B;M)$ and $(E;A,W;M)$ be double vector bundles with a
side bundle $A$ in common, and with cores $C$ and $L$ respectively. 
Suppose given a nondegenerate pairing 
$\lhangle\ ,\ \rhangle$ of $D$ over $A$ with $E$ over $A$, and two 
further nondegenerate pairings, both denoted $\langle\ ,\ \rangle$, 
of $B$ with $L$ and of $C$ with $W$, such that
\begin{enumerate}
\item for all $b\in B,\ \ell\in L$, 
$\lhangle \tilo^B_b,\Bar{\ell}\rhangle = \langle b, \ell\rangle$;
\item for all $c\in C,\ w\in W$, $\lhangle \Bar{c}, \tilo^W_w\rhangle
= \langle c, w\rangle$;
\item for all $c\in C,\ \ell\in L$, 
$\lhangle \Bar{c}, \Bar{\ell}\rhangle = 0;$
\item for all $d_1, d_2 \in D,\ e_1, e_2 \in E$ such that
$q^D_B(d_1) = q^D_B(d_2),$ \newline $q^E_W(e_1) = q^E_W(e_2),$
$q^D_A(d_1) = q^E_A(e_1),\ q^D_A(d_2) = q^E_A(e_2),$
we have 
$$
\lhangle d_1 \plusB d_2, e_1 \plusW e_2\rhangle =
\lhangle d_1, e_1 \rhangle + \lhangle d_2, e_2 \rhangle;
$$
\item for all $d\in D,\ e\in E$ such that $q^D_A(d) = q^E_A(e)$
and all $t\in\R$, we have
$$
\lhangle t\timesB d, t\timesW e\rhangle = t\lhangle d, e \rhangle.
$$
\end{enumerate}
(In all the above conditions we assume the various elements of the side
bundles lie in compatible fibres over $M$.)

Then the map $Z\co D\to E\duer A$ defined by
$\langle Z(d), e\rangle_A = \lhangle d,e\rhangle$ is an isomorphism 
of double vector bundles, with respect to $\id\co A\to A$ and the 
isomorphisms $B\to L^*$ and $C\to W^*$ induced by the pairings in 
\emph{(i)} and \emph{(ii)}.
\end{prop}

A pairing $\lhangle\ ,\ \rhangle$ satisfying the conditions of 
\ref{prop:dudvb} is called a \emph{pairing of the double vector 
bundles}.

Applying this result to the pairing (\ref{eq:3duals}) of 
$D\duer A$ and $D\duer B$, we find that the induced pairing of 
$B$ with $B^*$ is the standard one, but that of $A^*$ with $A$ is 
the negative of the standard pairing. Hence the signs in the 
following result are unavoidable. 

\begin{cor}                                
\label{cor:56}
The pairing \emph{(\ref{eq:3duals})} induces isomorphisms of double 
vector bundles
\begin{gather*}
Z_A\co D\duer A\to D\duer B\duer C^*,\qquad
\langle Z_A(\Phi), \Psi\rangle_{C^*} = \lhangle \Phi, \Psi\rhangle\\
Z_B\co D\duer B\to D\duer A\duer C^*,\qquad
\langle Z_B(\Psi), \Phi\rangle_{C^*} = \lhangle \Phi, \Psi\rhangle
\end{gather*}
with $(Z_A)\duer C^* = Z_B$. Both isomorphisms induce the identity 
on the sides $C^*\to C^*$.

$Z_A$ is the identity on the cores $B^*\to B^*$, and induces $-\id$
on the side bundles $A\to A$. 

$Z_B$ is the identity on the side bundles $B\to B$, and induces $-\id$
on the cores $A^*\to A^*$. 
\end{cor}

\begin{ex}\rm 
Consider a trivial double vector bundle 
$$
D = A\times_M B \times_M C.
$$
Let $\Phi = (a, \psi, \kappa)$ be an element of 
$D\duer A = A\times_M B^*\times_M C^*$ and let 
$\Psi = (\phi, b, \kappa)$ be an element of $D\duer B$. Then taking 
any $d = (a, b, c)\in D$, we find that
$$
\lhangle \Phi, \Psi\rhangle 
	= \langle\psi, b\rangle - \langle\phi, a\rangle. 
$$

The associated maps are given by
\begin{gather*}
Z_A\co A\pback{M} B^* \pback{M} C^* \to A\pback{M} B^* \pback{M} C^*,
 \quad 
(a,\psi,\kappa)\mapsto (-a, \psi, \kappa);\\
Z_B\co A^*\pback{M} B \pback{M} C^* \to A^*\pback{M} B \pback{M} C^*,
 \quad 
(\phi,b,\kappa)\mapsto (-\phi, b, \kappa).
\end{gather*}
\quend
\end{ex}

The following result is essentially equivalent to Theorem 
\ref{thm:dualduality}, but deserves independent statement. 

\begin{thm}
\label{thm:Q}
For any double vector bundle $(D;A,B;M)$ there is a canonical 
isomorphism $Q$ from $D$ to the flip of $(D\duer A\duer C^*\duer B)$ 
which preserves the side bundles $A$ and $B$ and is $-\id$ on the 
cores $C$. 
\end{thm}

\pf
Let $\Pi = Z_A\duer A$ be the dualization of $Z_A$ over $A$. 
Denote by $F\co D\to D$ the map $d\mapsto \minusB d$, and define
$Q = (F\circ \Pi)^{-1}.$
\ked

There are now three operations on double vector bundles: taking the
vertical dual, denoted by $V$, taking the horizontal dual, denoted
$H$, and the operation $VHV$ which by Theorem \ref{thm:Q} combines 
the flip and reversal of the sign on the core; we denote this by 
$P$. We have $V^2 = H^2 = P^2 = I$, the identity operation and, 
by the same method as Theorem \ref{thm:Q}, $HVH = P$. The group 
generated by $V$, $H$ and $P$ therefore has elements
\begin{equation}
\label{eq:dihedral}
I,\quad V,\quad HV,\quad VHV = P,\quad 
HVHV = HP = VH,\quad VHVHV = H, 
\end{equation}
and is the dihedral group $\Delta_3$ of the triangle, or
the symmetric group $\mathcal{S}_3$. 

\section{The duals of $TE$}                  
\label{sect:proldual}

Consider the tangent prolongation double vector bundle 
(Figure~\ref{fig:flip_TE}(b)) of a vector bundle $(E, q, M)$. 

First consider the horizontal dual. The canonical pairing of
$E^*$ with $E$ prolongs to a pairing of $T(E^*)\to TM$ with 
$TE\to TM$. Suppose given $\mathfrak{X}\in T(E^*)$ and 
$\xi\in TE$ with $T(q)(\mathfrak{X}) = T(q_*)(\xi)$. Then 
$\mathfrak{X} = \ddt{\phi_t}\in T(E^*)$ and $\xi = \ddt{e_t}\in TE$
where $e_t\in E$ and $\phi_t\in E^*$ can be taken so that 
$q_*(\phi_t) = q(e_t)$ 
for $t$ near zero. Now define the \emph{tangent pairing} 
$\llangle~,~\rrangle$  by
\begin{equation}                             
\label{eq:pairing}
\llangle\mathfrak{X},\xi\rrangle = \ddt{\langle\phi_t,e_t\rangle}.
\end{equation}
To show that this is non--degenerate it is sufficient to work 
locally. Suppose, therefore, that $E = M\times V$. 
Regard $\xi$ as $(x_0,v_0,w_0)\in T_{m_0}M\times V\times V$ and 
$\mathfrak{X}$ as $(x_0, \phi_0, \psi_0)\in 
T_{m_0}M\times V^*\times V^*$. Then
$$
\mathfrak{X} = \ddt{(m_t, \phi_0 + t\psi_0)},\qquad
\xi = \ddt{(m_t, v_0 + tw_0)},
$$
where $\ddt{m_t} = x = T(q)(\mathfrak{X}) = T(q_*)(\xi)$. So
$$
\llangle\mathfrak{X},\xi\rrangle = 
\ddt{\langle \phi_0 + t\psi_0, v_0 + tw_0\rangle}
$$
Expanding out the RHS, the constant term and the quadratic term 
vanish in the derivative, and we are left with
$$
\llangle\mathfrak{X},\xi\rrangle = 
\langle\psi_0, v_0\rangle + \langle\phi_0, w_0\rangle
$$
from which it is clear that $\llangle~,~\rrangle$ is 
non--degenerate. We now need to establish that this is a 
pairing of the double vector bundles. 

\begin{prop}
The tangent pairing $\llangle\ ,\ \rrangle$ of $T(E^*)$ with $TE$ 
over $TM$ satisfies the conditions of Proposition 
{\rm \ref{prop:dudvb}}. In particular, for $m\in M$ and 
$\phi, \phi_1,\phi_2\in E^*_m,\ e, e_1, e_2\in E_m$,
\begin{gather*}
\llangle\Bar{\phi},\Bar{e}\rrangle = 0,  \qquad 
\llangle\tilo_\phi,\tilo_e\rrangle = 0,\\
\llangle\tilo_\phi,\Bar{e}\rrangle = \langle\phi,e\rangle,  \qquad 
\llangle\Bar{\phi},\tilo_e\rrangle = \langle\phi,e\rangle.
\end{gather*}
\vspace{-8pt} 
and  
\begin{gather*}
\llangle\tau_*(\phi_1,\phi_2),\tau(e_1,e_2)\rrangle =
      \langle\phi_1, e_2\rangle + \langle\phi_2, e_1\rangle,
\end{gather*}
where $\tau_*$ and $\tau$ are the injections in the core sequences 
of $T(E^*)$ and $T(E)$.
\end{prop}

\pf
These are easily verified from the definition. For example, 
$\Bar{\phi} = \left.\frac{d}{dt}(t\phi)\right|_0$ and 
$\Bar{e} = \left.\frac{d}{dt}(te)\right|_0$ so 
$$
\llangle\Bar{\phi}, \Bar{e}\rrangle = 
\left.\frac{d}{dt} t^2\langle\phi, e\rangle\right|_0 = 0
$$
whereas $\tilo_\phi = \left.\frac{d}{dt} \phi\right|_0$ so
$$
\llangle\tilo_\phi, \Bar{e}\rrangle 
	= \left.\frac{d}{dt} t\langle\phi, e\rangle\right|_0 
	= \langle\phi, e\rangle.
$$
The bilinearity conditions are easily verified and the final equation
follows. 
\ked

Thus the pairing of the cores of $T(E^*)$ and $T(E)$ is the zero 
pairing, and so too is the pairing of the zero sections above $E^*$ 
and $E$. However the core of $T(E^*)$ and the zero section of 
$T(E)$ are paired under the standard pairing, and the same is true 
of the zero section of $T(E^*)$ and the core of $T(E)$. 

It now follows that there is an isomorphism of double vector bundles
from $T(E^*)$ to the dual $TE\duer TM$ of $TE$ over $TM$. For 
convenience we denote this simply by $T^\sol E$ and call it the 
\emph{prolongation dual of} $TE$. The next 
result follows from the general theory of \S\ref{sect:ddvb}. 

\begin{prop}{\rm \citep{MackenzieX:1994}}
\label{prop:I}
The map $I\co T(E^*)\to T^\sol(E)$ defined by
$$
\langle I(\mathcal{X}), \eta\rangle_{TM} 
	= \llangle \mathcal{X}, \eta\rrangle
$$
where $\mathcal{X}\in T(E^*),\ \eta\in TE$, is an isomorphism of 
double vector bundles preserving the side bndles $E^*$ and $TM$ and 
the core bundles $E^*$. 
\end{prop}

When a name is needed we call $I$ the 
\emph{internalization map}. 
In future we will almost always work with $T(E^*)$ and the tangent 
pairing, rather than with $T^\sol E$ and $I$. 
 
Now consider the vertical dual of $TE$. Since the core of the 
double vector bundle $TE$ is $E$, dualizing the 
structure over $E$ leads to a double vector bundle of the form
\begin{equation}                         \label{eq:cotdual}
\begin{matrix}
        &&r_E&&\\
        &T^*E&\mlra&E^*&\\
        &&&&\\
     c_E&\Big\downarrow&&\Big\downarrow&q_{E^*}\\
        &&&&\\
        &E&\mlra&M;&\\
        &&q_E&&\\
\end{matrix}
\end{equation}
We refer to this as the \emph{cotangent dual} of $TE$. We will 
give a detailed description of the structures involved. Although 
this is a special case of the general construction, this example is 
so basic to the rest of the paper that it merits a specific 
treatment. 

In (\ref{eq:cotdual}) the vertical bundle is the standard cotangent 
bundle
of $E$, and the notation $T_X^*(E)$ will always refer to the fibre 
with respect to $c_E$. In this bundle we use standard notation, and 
denote the zero element of $T_X^*(E)$ by $\tilo^*_X$. We drop the 
subscripts $E$ from the maps when no confusion is likely.

The map $r\co T^*E\to E^*$ takes the form
\begin{equation}
\label{eq:rdf}
\langle r(\Phi),Y\rangle 
= \langle\Phi,\tau(X,Y)\rangle 
= \langle\Phi, \tilo_X\dpl\Bar{Y}\rangle
\end{equation}
where $\Phi\in T_X^*(E)$,  $X\in E_m$ and $Y\in E_m$. Thus 
$r(\Phi)\in E^*_m$. For the addition over $E^*$ we have
$$
\langle\Phi\plusEst\Psi, \xi\dpl\eta\rangle 
= \langle\Phi, \xi\rangle + \langle\Psi,\eta\rangle,
$$
where
$\Phi\in T_X^*(E),\ \Psi\in T_Y^*(E)$ with $r(\Phi) = r(\Psi)\in
E^*_m$, and $\xi\in T_X(E),\ \eta\in T_Y(E)$ with
$T(q)(\xi) = T(q)(\eta)$.
This defines $\Phi\plusEst\Psi\in T_{X+Y}^*(E)$. 
Similarly, we have
$$
\langle t\timesEst\Phi, t\dtimes \xi\rangle 
= t\langle\Phi,\xi\rangle,
$$
for $t\in\R$ and $\xi\in T_X(E)$. The zero element of 
$r^{-1}(\phi)$, where $\phi\in E^*_m$, is 
$\tilo^r_\phi\in T_{0_m}^*(E)$ where
$$
\langle\tilo^r_\phi, T(0)(x) + \Bar{X}\rangle 
= \langle\phi,X\rangle
$$
for $x\in T_m(M), X\in E_m$.

Given $\omega\in T_m^*(M)$, the corresponding core element 
$\Bar{\omega}$ is
$$
\langle\Bar{\omega}, T(0)(x) + \Bar{X}\rangle 
= \langle\omega,x\rangle,
$$
for $x\in T_m(M),\ X\in E_m.$ The injection over $E$,
$$
q\xclam T^*M\to T^*E,\ 
(X,\omega)\mapsto \tilo^*_X\plusEst\Bar{\omega},
$$
is precisely the dual of $T(q)\xclam$; that is to say, it is the map 
corresponding to the lifting of 1--forms from $M$ to $E$. Thus 
$\tilo^*_X\plusEst\Bar{\omega}$ is the pullback of 
$\omega\in T_m^*(M)$ to $E$ at the point $X\in E_m$. 

The core exact sequence for $c$ is
\begin{equation}
q\xclam T^*M\inj T^*E \stackrel{r\xclam = \tau^*}{\surj} q\xclam E^*,
\end{equation}
and this is the dual of the core exact sequence for $TE$ and $p_E$. 
The other core exact sequence is
\begin{equation}
q_*\xclam T^*M\inj T^*E \stackrel{c\xclam}{\surj} q_*\xclam E,
\end{equation}
where each bundle here is over $E^*$. The injection
$q_*\xclam T^*M\to T^*E$
is $(\phi,\omega)\mapsto\tilo^r_\phi + \Bar{\omega}$ and
$\langle\tilo^r_\phi + \Bar{\omega},T(0)(x) + \Bar{X}\rangle 
= \langle\phi, X\rangle + \langle\omega,x\rangle.$ 

Given $\omega\in T^*M$, the corresponding core element is
$$
\Bar{\omega} = (\omega,0,0).
$$

To summarize, the two dual double vector bundles of $D = TE$ are
$$     %
\begin{matrix}
        &&      &&\\
D\duer E =&T^*E&\longrightarrow&E^*&\\
        &&&&\\
        &\Big\downarrow&&\Big\downarrow&\\
        &&&&\\
        &E&\longrightarrow&M&\\
        &&&&\\
\end{matrix}
\mbox{\quad and\quad}
\begin{matrix}
        &&      &&\\
D\duer TM =&T^\sol E&\longrightarrow&TM&\\
        &&&&\\
        &\Big\downarrow&&\Big\downarrow&\\
        &&&&\\
        &E^*&\longrightarrow&M&\\
        &&&&\\
\end{matrix}
$$                
and the pairing
\begin{equation}
\label{eq:pair}
\lhangle\Phi, \gof\rhangle = 
\langle\Phi, \xi\rangle_E - \langle\gof, \xi\rangle_{TM}
\end{equation}
for suitable $\xi\in TE$. Composing the isomorphism $Z_E$ from 
\ref{cor:56} with the dual over $E^*$ of the internalization
isomorphism $I$, we get an isomorphism of double vector bundles
$$
(I\duer E^*)\circ Z_E\co T^*E\to T^*(E^*);
$$
denote this temporarily by $S^{-1}$. For $\Phi\in T^*E$ we have
\begin{multline*}
\langle S^{-1}(\Phi), \mathcal{X}\rangle_{E^*} = 
\langle (I\duer E^*)\circ Z_E(\Phi), \mathcal{X}\rangle_{E^*} = 
\langle Z_E(\Phi), I(\mathcal{X})\rangle_{E^*}\\ = 
\lhangle \Phi, I(\mathcal{X})\rhangle = 
\langle\Phi, \xi\rangle_E - \langle I(\mathcal{X}), \xi\rangle_{TM} 
= \langle \Phi, \xi\rangle_E - \llangle \mathcal{X}, \xi\rrangle.
\end{multline*}
Here we used the definition of $Z_E$, the definition 
(\ref{eq:pair}), and the definition of $I$. It follows that 
for $\goF\in T^*(E^*)$, writing $\Phi = S(\goF)$, we have 
$$
\langle \goF, \mathcal{X}\rangle_{E^*} = 
\langle S(\goF), \xi\rangle_E - \llangle \mathcal{X}, \xi\rrangle.
$$
Recall that $I$, and hence its dual, preserves both sides and the 
core, whereas $Z_E$ induces $-\id$ on the sides $E$. We therefore 
define
$$
R\co T^*(E^*)\to T^*(E),\qquad R(\goF) = S(\minusEst\goF). 
$$
To summarize:

\begin{thm}{\rm \citep{MackenzieX:1994}}
\label{thm:MX55}
The map $R$ just defined is an isomorphism of double vector 
bundles, preserving the side bundles $E$ and $E^*$, and inducing 
$-\id\co T^*M\to T^*M$ on the cores. Further, for all 
$\xi\in TE,\ \mathcal{X}\in T(E^*),\ \goF\in T^*(E^*)$ such that 
$\xi$ and $\mathcal{X}$ have the same projection into $TM$, 
$\mathcal{X}$ and $\goF$ have the same projection into $E^*$, 
and $\goF$ and $\xi$ have the same projection into $E$, 
\begin{equation}
\label{eq:MX38}
\llangle \mathcal{X}, \xi\rrangle = 
	\langle R(\goF), \xi\rangle_E
		+ \langle \goF, \mathcal{X}\rangle_{E^*}  
\end{equation}
\end{thm}

We call $R$ the \emph{reversal isomorphism}. It is proved in
\citep{MackenzieX:1994} that $R$ is an antisymplectomorphism of
the canonical symplectic structures.  

\section{Triple vector bundles}
\label{sect:tvbs}

The definition of a triple vector bundle follows the same pattern
as in the double case. There are a number of evident 
reformulations. 

\begin{df}
A \emph{triple vector bundle} is a manifold $\tvb$ together with
three vector bundle structures, over bases $D_1, D_2, D_3$, each
of which is a double vector bundle with side bundles respectively
$E_2$ and $E_3$, $E_3$ and $E_1$, $E_1$ and $E_2$, where
$E_1, E_2, E_3$ are vector bundles over a shared base $M$, such that
each pair of vector bundle structures on $\tvb$ forms a double
vector bundle, the operations of which are vector bundle 
morphisms with respect to the third vector bundle structure. 
\end{df}

We display a triple vector bundle in a diagram such as 
Figure~\ref{fig:tvb}(a). (We always read figures of this type with
$(\tvb;D_1, D_2;E_3)$ at the rear and $(D_3;E_2,E_1;M)$ coming 
out of the page toward the reader.) The three structures of 
double vector bundle on $\tvb$ are the \emph{upper} double vector 
bundles, and $D_1, D_2, D_3$ are the \emph{lower} double vector 
bundles. We refer to $(D_1;E_2, E_3;M)$ as the \emph{floor} of
$\tvb$ and to $(\tvb;D_2, D_3;E_1)$ as the \emph{roof} of $\tvb$. 

\begin{figure}[htb]
\setlength{\unitlength}{1cm}
\begin{picture}(5,4)
\put(1,3){$\begin{matrix}
                      &&                        &\\
                      &\tvb&\longrightarrow &D_2\\
                      &&&\\
                      &\Big\downarrow& &\Big\downarrow\\
                      &&&\\
                      &D_1&\longrightarrow &E_3\\
            \end{matrix}
                      $}
\put(1.8, 3.5){\vector(2,-1){1}}                    
\put(2.9, 3.5){\vector(2,-1){1}}                    
\put(1.8, 1.8){\vector(2,-1){1}}                    
\put(2.9, 1.8){\vector(2,-1){1}}                    

\put(2.7,1.7){$\begin{matrix}
                     &&                        &\\
                     &D_3&\longrightarrow &E_1\\
                     &&&\\
                     &\Big\downarrow&&\Big\downarrow\\
                     &&&\\
                     &E_2&\longrightarrow & M\\
		     &&&\\
		     &\mbox{(a)}&&\\	
               \end{matrix}
                     $}

\put(7,3){$\begin{matrix}
                      &&                        &\\
                      &\Bar{k}_1&\longrightarrow &d_2\\
                      &&&\\
                      &\Big\downarrow& &\Big\downarrow\\
                      &&&\\
                      &d_1&\longrightarrow &e_3\\
            \end{matrix}
                      $}
\put(7.8, 3.5){\vector(2,-1){1}}                    
\put(8.9, 3.5){\vector(2,-1){1}}                    
\put(7.8, 1.8){\vector(2,-1){1}}                    
\put(8.9, 1.8){\vector(2,-1){1}}                    
\put(8.7,1.7){$\begin{matrix}
                     &&                        &\\
                     &d_3&\longrightarrow &e_1\\
                     &&&\\
                     &\Big\downarrow&&\Big\downarrow\\
                     &&&\\
                     &e_2&\longrightarrow & m\\
		     &&&\\
		     &\mbox{(b)}&&\\	
               \end{matrix}
                     $}
\end{picture}
\caption{\ \label{fig:tvb}} 
\end{figure}

We have found that, rather than assembling a notation capable of 
handling any calculation in a triple vector bundle without 
ambiguity, 
it is generally preferable to develop an ad hoc notation for each
occasion. The great majority of calculations use only certain parts
of the structure, and in such cases a modification of the 
notation of \S\ref{sect:dvbs} is often sufficient. 

Each of the lower double vector bundles $D_i$ has a core, which 
is denoted $C_i$. The core of the upper double vector bundle 
$(\tvb;D_3, D_2; E_1)$ is denoted $K_1$. Consider a core element
$\Bar{k}_1\in\tvb$ where $k_1$ projects to $e_1\in E^1_m$. Then 
the $d_2$
in Figure~\ref{fig:tvb}(b) %
is the zero over $e_1$ for $D_2\to E_1$, and $d_3$ is the
zero over $e_1$ for $D_3\to E_1$. From the morphism condition we
then have that $e_2$ and $e_3$ are zeros over $m$. So 
$d_1 = \Bar{c}$ is a core element for some $c\in C_1$. This defines
a map $K_1\to C_1$. For $k, k'\in K_1$ over the same element of 
$C_1$, define
\begin{equation}
\label{eq:crossadd}
\Bar{k +_{C_1} k'} = \Bar{k} +_{D_1} \Bar{k'}, 
\end{equation}
where each of the three bars refers to the roof double vector 
bundle. With scalar multiplication defined in a similar fashion, 
$K_1\to C_1$ is a vector bundle, and a double vector
bundle as shown in Figure~\ref{fig:cdvbs}(a). 

The cores of the other 
upper double vector bundles are likewise denoted $K_2$ and $K_3$ and 
form double vector bundles as in Figure~\ref{fig:cdvbs}(b)(c).
These three are the \emph{core double vector bundles}. Although 
defined by restrictions of the operations in $\tvb$, they are not 
substructures of $\tvb$. 

\begin{figure}[htb]
\setlength{\unitlength}{1cm}
\begin{picture}(10,3.4)
\put(1,1.6){$\begin{matrix}
                      &&                        &\\
                      &K_1&\longrightarrow &E_1\\
                      &&&\\
                      &\Big\downarrow& &\Big\downarrow\\
                      &&&\\
                      &C_1&\longrightarrow &M\\
                      &&&\\
		      &&\mbox{(a)}&\\
            \end{matrix}
                      $}
\put(5,1.61){$\begin{matrix}
                      &&                        &\\
                      &K_2&\longrightarrow &E_2\\
                      &&&\\
                      &\Big\downarrow& &\Big\downarrow\\
                      &&&\\
                      &C_2&\longrightarrow &M\\
                      &&&\\
		      &&\mbox{(b)}&\\
            \end{matrix}
                      $}
\put(9,1.6){$\begin{matrix}
                      &&                        &\\
                      &K_3&\longrightarrow &E_3\\
                      &&&\\
                      &\Big\downarrow& &\Big\downarrow\\
                      &&&\\
                      &C_3&\longrightarrow &M\\
                      &&&\\
		      &&\mbox{(c)}&\\
            \end{matrix}
                      $}
\end{picture}
\caption{\ \label{fig:cdvbs}}
\end{figure}

Denote the core of $(K_1;C_1, E_1;M)$ by $W$. In 
Figure~\ref{fig:tvb}(b), let $k_1 = \Bar{w}$ where $w\in W$. Then
$c_1$ is the zero of $C_1$ over $m$, and so $d_1 = \Bar{c}_1$ is
a double zero of $D_1$. 

Next, since $e_2$ is a zero, $(d_3;e_2, e_1;m)$ must be of the
form $d_3 = 0^{D_3}_{e_1} +_{E_2} \Bar{c}_3$. But $e_1$ is zero, 
since $w$ is in the core, so $d_3$ is a core element. But it is
known to be a zero over $e_1$, so must be a double zero. Similarly
$d_2$ is a double zero. 

This proves most of the next result, and the remainder is an
easy verification. 

\begin{prop}
\label{prop:threecores}
Each of the core double vector bundles has as core the set
$W$ of elements $w\in\tvb$ for which the projection to each of
$D_1, D_2, D_3$ is a double zero of the lower double vector bundle. 
Further, the vector bundle structures with base $M$ induced 
on $W$ by the various core double vector bundles coincide. 

Given $w\in W$, the core elements in $\tvb$ corresponding to
$\Bar{w}^1\in K_1$, $\Bar{w}^2\in K_2$ and $\Bar{w}^3\in K_3$,
coincide. 
\end{prop}

We call $W$ the \emph{ultracore of} $\tvb$. 

\begin{ex}\rm
For a double vector bundle $(D;A,B;M)$, the \emph{tangent 
prolongation triple vector bundle} is as shown in 
Figure~\ref{fig:cotvb}(a). 
Two of the three core double vector bundles are copies of
$(D;A,B;M)$ and the third is $(TC;C,TM;M)$. The ultracore 
of $TD$ is $C$, the core of $D$. It is illuminating to
verify (\ref{eq:crossadd}) and \ref{prop:threecores} directly
in this example. 
\quend\end{ex}

\begin{figure}[htb]
\setlength{\unitlength}{1cm}
\begin{picture}(5,4)
\put(1,3){$\begin{matrix}
                      &&                        &\\
                      &TD&\longrightarrow &TB\\
                      &&&\\
                      &\Big\downarrow& &\Big\downarrow\\
                      &&&\\
                      &D&\longrightarrow &B\\
            \end{matrix}
                      $}
\put(1.8, 3.5){\vector(2,-1){1}}                    
\put(3, 3.5){\vector(2,-1){1}}                    
\put(1.8, 1.8){\vector(2,-1){1}}                    
\put(3, 1.8){\vector(2,-1){1}}                    

\put(2.7,1.7){$\begin{matrix}
                     &&                        &\\
                     &TA&\longrightarrow &TM\\
                     &&&\\
                     &\Big\downarrow&&\Big\downarrow\\
                     &&&\\
                     &A&\longrightarrow & M\\
		     &&&\\
		     &\mbox{(a)}&&\\	
               \end{matrix}
                     $}

\put(7,3){$\begin{matrix}
                      &&                        &\\
                      &T^*D&\longrightarrow &D\duer B\\
                      &&&\\
                      &\Big\downarrow& &\Big\downarrow\\
                      &&&\\
                      &D&\longrightarrow &B\\
            \end{matrix}
                      $}
\put(7.8, 3.5){\vector(2,-1){1}}                    
\put(7.8, 1.8){\vector(2,-1){1}}                    
\put(9.1, 1.8){\vector(2,-1){1}}                    
\put(8.6,1.7){$\begin{matrix}
                     &&                        &\\
                     &\phantom{X}D\duer A&&\\
                     &&&\\
                     &\Big\downarrow&&\\
                     &&&\\
                     &A&\longrightarrow & M\\
		     &&&\\
		     &\mbox{(b)}&&\\	
               \end{matrix}
                     $}
\end{picture}
\caption{\ \label{fig:cotvb}}
\end{figure}

The key to understanding the relations between the duals of a
double vector bundle, and the role of the dihedral group
(\ref{eq:dihedral}), lies in constructing a cotangent form of
this example. Both the left and the rear faces of 
Figure~\ref{fig:cotvb}(a) are tangent prolongation double vector 
bundles of ordinary vector bundles and so we can form the figure 
in Figure~\ref{fig:cotvb}(b). 

In Figure~\ref{fig:cotstage1}(a) we have added the two double vector
bundle duals associated to $D$. In this diagram each of the four
vertical sides, and of course the floor, is a double vector
bundle. We need to prove that the roof is a double vector bundle
and that Figure~\ref{fig:cotstage1}(a) is a triple vector bundle. 

First consider the roof, shown in Figure~\ref{fig:cotstage1}(b). 
We use a short notation for the projections. 

\begin{figure}[htb]
\setlength{\unitlength}{1cm}
\begin{picture}(5,4)
\put(1,3){$\begin{matrix}
                      &&                        &\\
                      &T^*D&\longrightarrow &D\duer B\\
                      &&&\\
                      &\Big\downarrow& &\Big\downarrow\\
                      &&&\\
                      &D&\longrightarrow &B\\
            \end{matrix}
                      $}
\put(1.8, 3.5){\vector(2,-1){1.4}}                    
\put(3.1, 3.5){\vector(2,-1){1.4}}                    
\put(1.8, 1.8){\vector(2,-1){1.4}}                    
\put(3.1, 1.8){\vector(2,-1){1.4}}                    

\put(3,1.4){$\begin{matrix}
                     &&                        &\\
                     &D\duer A&\longrightarrow &C^*\\
                     &&&\\
                     &\Big\downarrow&&\Big\downarrow\\
                     &&&\\
                     &A&\longrightarrow & M\\
		     &&&\\
		     &\mbox{(a)}&&\\	
               \end{matrix}
                     $}
\put(8,2){$
\begin{matrix}
 &&    r_B             &&\\
 &T^*D&\longrightarrow &D\duer B&\\
 &&&&\\
r_A &\Big\downarrow& &\Big\downarrow& q_{\du B}\\
 &&&&\\
 &D\duer A&\longrightarrow &C^*&\\
 && q_{\du A}        &&\\
 &&&&\\
 &&\mbox{(b)}&&\\
\end{matrix}
$}
\end{picture}
\caption{\ \label{fig:cotstage1}}  
\end{figure}

In this structure, we know that each side is a vector bundle. 

\begin{prop}
\label{prop:roofisdvb}
The structure in Figure {\rm \ref{fig:cotstage1}(b)} is a double 
vector bundle with core $T^*C$. 
\end{prop}

\pf
First we must prove that the projections form a commutative
square. Take $\gof\in T_d^*D$, where $d$ has the form
$(d;a,b;m)$. Then, for all $(d';a, b'; m)$, 
$$
\langle r_A(\gof), d'\rangle 
= \langle\gof, \tilo_d \DPL{A} \Bar{d'}^A\rangle.
$$
Here $\tilo_d$ is the zero of $TD\to D$ above $d$ and 
the subscript on $\DPL{A}$ indicates that this is the tangent of
the addition in $D\to A$. The superscript $A$ on the bar 
indicates that $D$ is here the core of $(TD;TA,TM;A).$

Writing $\phi = r_A(\gof)$ we next have
$$
\langle q_{\du A}(\phi), c\rangle 
= \langle\phi, 0_a^D +_B \Bar{c}\rangle. 
$$
Here $0_a^D$ is the zero of $D\to A$ over $a$ and 
$\Bar{c}$ is the core element of $D$ corresponding to $c$. 
Writing out the corresponding formulas for the other side, we
must prove that
\begin{equation}
\label{eq:commute}
\tilo_d \DPL{A} \Bar{(0_a^D+_B\Bar{c})}^A
= \tilo_d \DPL{B} \Bar{(0_b^D+_A\Bar{c})}^B.
\end{equation}
Using (\ref{eq:crossadd}), 
the LHS becomes
$$
\tilo_d \DPL{A} (\Bar{0_a^D}^A \DPL{B} \Bar{\Bar{c}}\,)
$$
Writing in terms of tangents to curves, we have
$$
\tilo_d = \ddt{d},\qquad
\Bar{0^D_a}^A = \ddt{(t\dtimes_A 0^D_a)} = \ddt{0^D_a},\qquad
\Bar{\Bar{c}} = \ddt{t\,\Bar{c}}.
$$
Now, using the interchange rule, 
\begin{multline*}
d +_A (0^D_a +_B \Bar{c}) 
= (d +_B 0^D_b) +_A (0^D_a +_B \Bar{c}) \\
= (d +_A 0^D_a) +_B (0^D_b +_A \Bar{c}) 
= d +_B (0^D_b +_A \Bar{c}) 
\end{multline*}
and from this (\ref{eq:commute}) follows. The proof that 
$r_A$ preserves the addition and scalar multiplication proceeds 
in a similar way. 

Next we show that the core is $T^*C$. Suppose that $\gof\in T^*_dD$
maps to zero under both $r_A$ and $r_B$. Then $d = \Bar{c}$ is a
core element and $\gof$ vanishes on elements $\xi\in T_{\Bar{c}}D$
which are vertical with respect to either $q^D_A$ or $q^D_B$. 
If $\xi$ is vertical with respect to $q^D_A$ then, in the notation
of Figure~\ref{fig:proof}(a), $X = \tilo^A_m$, and it follows 
that $Z = \tilo_m$ and so $Y$ is a core element. Likewise if
$\xi$ is vertical with respect to $q^D_B$, then $Y = \tilo^B_m$
and $X$ is a core element. Adding two such representative elements,
it follows that $\gof$ vanishes on all $\mathcal{X}$ as shown in
Figure~\ref{fig:proof}(b),

Now take $\xi\in T_{\Bar{c}}D$ as shown in 
Figure~\ref{fig:proof}(a). 
\begin{figure}[htb]
\setlength{\unitlength}{1cm}
\begin{picture}(5,4.5)
\put(1,3){$\begin{matrix}
                      &&                        &\\
                      &\xi&\longrightarrow &Y\\
                      &&&\\
                      &\Big\downarrow& &\Big\downarrow\\
                      &&&\\
                      &\Bar{c}&\longrightarrow &0^B_m\\
            \end{matrix}
                      $}
\put(1.5, 3.5){\vector(2,-1){1}}                    
\put(2.7, 3.5){\vector(2,-1){1}}                    
\put(1.5, 1.8){\vector(2,-1){1}}                    
\put(2.7, 1.8){\vector(2,-1){1}}                    

\put(2.4,1.7){$\begin{matrix}
                     &&                        &\\
                     &X&\longrightarrow &Z\\
                     &&&\\
                     &\Big\downarrow&&\Big\downarrow\\
                     &&&\\
                     &0^A_m&\longrightarrow & m\\
		     &&&\\
		     &\mbox{(a)}&&\\	
               \end{matrix}
                     $}

\put(7,3){$\begin{matrix}
                      &&                        &\\
                      &\mathcal{X}&\longrightarrow &\Bar{b}\\
                      &&&\\
                      &\Big\downarrow& &\Big\downarrow\\
                      &&&\\
                      &\Bar{c}&\longrightarrow &0^B_m\\
            \end{matrix}
                      $}
\put(7.5, 3.5){\vector(2,-1){1}}                    
\put(8.5, 3.5){\vector(2,-1){1}}                    
\put(7.5, 1.8){\vector(2,-1){1}}                    
\put(8.5, 1.8){\vector(2,-1){1}}                    
\put(8.4,1.7){$\begin{matrix}
                     &&                        &\\
                     &\Bar{a}&\longrightarrow &\tilo_m\\
                     &&&\\
                     &\Big\downarrow&&\Big\downarrow\\
                     &&&\\
                     &0^A_m&\longrightarrow & m\\
		     &&&\\
		     &\mbox{(b)}&&\\	
               \end{matrix}
                     $}
\end{picture}
\caption{\ \label{fig:proof}} 
\end{figure}
Because $X\in T_{0^A_m}A$, it has the form $X = T(0^A)(Z) + \Bar{a}$
for some $a\in A_m$; likewise $Y$ has the form
$Y = T(0^B)(Z) + \Bar{b}$ for some $b\in B_m$. Now define
$$
\mathcal{X} = (T(0^D_B)(\Bar{b})\DPL{A}\tilo_{\Bar{c}})
               \underset{D}{+}
              (T(0^D_A)(\Bar{a})\DPL{B}\tilo_{\Bar{c}}).
$$
Then $\Bar{\gamma} = 
\xi \underset{D}{-}\mathcal{X}$ has the form shown in
Figure~\ref{fig:proof2}(a) and is an element of $\Bar{TC}$. 

It is now possible to extend a given $\omega\in T^*_cC$ 
to $\Bar{\omega}\in T^*_{\Bar{c}}D$ by
$$
\langle\Bar{\omega}, \xi\rangle 
= \langle\omega, \xi \underset{D}{-}\mathcal{X}\rangle
$$
and $\Bar{\omega}$ is annulled by both $r_A$ and $r_B$. This
$\Bar{\omega}$ is the \emph{core element} of the double vector
bundle in Figure~\ref{fig:cotstage1}(b). 

Now that the commutativity of the projections has been established, 
verification that addition and scalar multiplication are preserved
is straightforward. 
\ked

The double vector bundle in Figure~\ref{fig:cotstage1}(b) is
of a type not previously encountered. 

\begin{figure}[htb]
\setlength{\unitlength}{1cm}
\begin{picture}(5,4.5)
\put(1,3){$\begin{matrix}
                      &&                        &\\
                      &\Bar{\gamma}&\longrightarrow &T(0^B)(Z)\\
                      &&&\\
                      &\Big\downarrow& &\Big\downarrow\\
                      &&&\\
                      &\Bar{c}&\longrightarrow &0^B_m\\
            \end{matrix}
                      $}
\put(1.5, 3.5){\vector(2,-1){1.7}}                    
\put(3, 3.5){\vector(2,-1){1.7}}                    
\put(1.5, 1.8){\vector(2,-1){1.7}}                    
\put(3, 1.8){\vector(2,-1){1.7}}                    

\put(2.9,1.2){$\begin{matrix}
                     &&                        &\\
                     &T(0^A)(Z)&\longrightarrow &Z\\
                     &&&\\
                     &\Big\downarrow&&\Big\downarrow\\
                     &&&\\
                     &0^A_m&\longrightarrow & m\\
		     &&&\\
		     &\mbox{(a)}&&\\	
               \end{matrix}
                     $}

\put(7,3){$\begin{matrix}
                      &&                        &\\
                      &&&D\duer B\\
                      &&&\\
                      && &\Big\downarrow\\
                      &&&\\
                      &D&\longrightarrow &B\\
            \end{matrix}
                      $}
\put(9, 3.5){\vector(2,-1){1}}                    
\put(7.6, 1.8){\vector(2,-1){1}}                    
\put(8.9, 1.8){\vector(2,-1){1}}                    
\put(8.2,1.7){$\begin{matrix}
                     &&                        &\\
                     &\phantom{X}D\duer A&\longrightarrow &C^*\\
                     &&&\\
                     &\Big\downarrow&&\Big\downarrow\\
                     &&&\\
                     &A&\longrightarrow & M\\
		     &&&\\
		     &\mbox{(b)}&&\\	
               \end{matrix}
                     $}
\end{picture}
\caption{\ \label{fig:proof2}} 
\end{figure}

It is now a straightforward matter to complete the proof 
of the following:

\begin{thm}
\label{thm:T*D}
The structure in Figure~{\rm \ref{fig:cotstage1}(a)} is a triple 
vector bundle. 
\end{thm}

The core double vector bundles are $(T^*A;A, A^*;M)$, 
$(T^*B;B,B^*;M)$ and $(T^*C;C,C^*;M)$, and the ultracore is
$T^*M$. Observe that the triple vector bundle $T^*D$ has a much 
higher degree of symmetry than $TD$. 

\begin{ex}\rm
Consider seven
vector bundles $E_1, E_2, E_3, C_1, C_2, C_3, W$ over
a shared base $M$. Let $D_1$ be the trivial double vector bundle
with sides $E_2$ and $E_3$ and core $C_1$, and likewise form 
$D_2$ and
$D_3$. Similarly, let $K_1$ be the trivial double vector bundle
with sides $C_1$ and $E_1$ and core $W$, and form $K_2$ and $K_3$
in the same way. Lastly, let $\tvb$ be the pullback of all seven
vector bundles over $M$. Then $\tvb$ can be considered as the
trivial double vector bundle with side bundles $D_1\to E_3$ and
$D_2\to E_3$ and core $K_3\to E_3$. Likewise, $\tvb$ can be
considered the trivial double vector bundle over $D_2$ and $D_3$
with core $K_1$ and over $D_3$ and $D_1$ with core $K_2$. With
these structures, $\tvb$ is a triple vector bundle, the 
\emph{trivial triple vector bundle} determined by the given seven 
vector bundles. 

More refined versions of this construction exist. For example, 
suppose given four vector bundles $E_1, E_2, E_3, W$ on $M$ and 
three double vector bundles $(D_1;E_2, E_3;M)$, $(D_2;E_3,E_1;M)$, 
$(D_3;E_1,E_2;M)$. Then 
there is a triple vector bundle for which $D_1$, $D_2$, $D_3$ are
the lower double vector bundles and $W$ is the ultracore, and for 
which each of the core double vector bundles is trivial. 
\quend
\end{ex}

\section{Cornerings}
\label{sect:corners}

Continue with a double vector bundle $D$ as in the previous section. 
Since $D$ is a vector bundle over $A$, we have 
$T^*D\isom T^*(D\duer A)$ by \ref{thm:MX55}, and similarly 
$T^*D\isom T^*(D\duer B)$. Once it has been shown that
these isomorphisms respect the triple structures, we can regard
$\Delta_3$ as acting on the cube $T^*D$ by rotations about the
axis from $T^*D$ to $M$. 

\begin{thm}
The map $R^{-1}\co T^*D\to T^*(D\duer A)$ arising from the vector 
bundle $D\to A$ is an isomorphism of 
triple vector bundles over $Z_B\co D\duer B\to (D\duer A\duer C^*)$, 
the other maps on the side structures being identities. 
\end{thm}

\pf
The main work is to show that $R^{-1}$ is a morphism of
vector bundles over $Z_B$. Take $\Phi\in T^*_dD$ and denote
$R^{-1}(\Phi)$ by $\gof$. Let $d$ have the form $(d;a,b,m)$ and 
let the projections of $\Phi$ to $D\duer A$, $D\duer B$ and $C^*$ 
be $\chigh$, $\psi$ and $\kappa$ respectively. Since $R$ preserves
$D$ and $D\duer A$, it follows that $\gof$ projects to $d\in D$ and 
to $\chigh\in D\duer A$. 

For $\psi$ we have, from (\ref{eq:rdf}), 
\begin{equation}
\label{eq:rz1}
\langle\psi, d_1\rangle_B 
= \langle\Phi, \tilo_d\DPL{B}\Bar{d}_1^B\rangle_D
\end{equation}
for any $d_1$ of the form $(d_1;a_1, b;m)$. For $\gof$ and $\Phi$
we have, by (\ref{eq:MX38}), 
\begin{equation}
\label{eq:rz2}
\llangle\mathcal{X}, \xi\rrangle_{TA}
= \langle\Phi, \xi\rangle_D 
       + \langle\gof, \mathcal{X}\rangle_{D\duer A}
\end{equation}
where $\mathcal{X}\in T(D\duer A)$ has the form
$(\mathcal{X}; \chigh, X;a)$ for some $X\in TA$, and $\xi\in TD$ then
has the form $(\xi;d,X;a)$. Next, for $Z_B(\psi)\in
D\duer A\duer C^*$, we have, for each $\phi\in D\duer A$ of the 
form $(\phi;a_2,\kappa;m)$, 
\begin{equation}
\label{eq:rz3}
\langle Z_B(\psi), \phi\rangle = \lhangle\phi,\psi\rhangle_{C^*}
= \langle\phi, d_2\rangle_A - \langle\psi, d_2\rangle_B
\end{equation}
for any $d_2\in D$ of the form $(d_2;a_2, b; m)$. Lastly, for the
same $\phi$ we have
\begin{equation}
\label{eq:rz4}
\langle r_{C^*}(\gof), \phi\rangle 
= \langle\gof, 
\tilo^{(D\duer A)}_\chigh \DPL{C^*}\Bar{\phi}^{(D\duer A)}\rangle
\end{equation}
Here $\tilo^{(D\duer A)}_\chigh$ is the zero in $T(D\duer A)$ over
$\chigh$ and $\Bar{\phi}^{(D\duer A)}$ is the core element of
$T(D\duer A)$ corresponding to $\phi$. The addition is in the
bundle $T(D\duer A)\to TC^*$. 

We must prove that the RHSs of (\ref{eq:rz3}) and (\ref{eq:rz4})
are equal. We substitute  
$$
\mathcal{X} 
= \tilo^{(D\duer A)}_\chigh\DPL{C^*}\Bar{\phi}^{(D\duer A)}
\quad\mbox{and}\quad
\xi = \tilo_d\DPL{B}\Bar{d}_1^B
$$
into (\ref{eq:rz2}). Providing $a_1 = a_2$, the relevant 
projections match. Now applying (iv) of \ref{prop:dudvb} 
to the double vector bundles $(T(D\duer A); TA, TC^*; TM)$ 
and $(TD;TA, TB;TM)$, we have
\begin{multline*}
\llangle\mathcal{X}, \xi\rrangle_{TA}
= \llangle 
  \tilo^{(D\duer A)}_\chigh\DPL{C^*}\Bar{\phi}^{(D\duer A)},\ 
  \tilo_d\DPL{B}\Bar{d}_1^B
  \rrangle_{TA}\\
= \llangle\tilo^{(D\duer A)}_\chigh,\ \tilo_d \rrangle_{TA}
  + \llangle \Bar{\phi}^{(D\duer A)},\Bar{d}_1^B \rrangle_{TA}.
\end{multline*}
In the first term, $\tilo^{(D\duer A)}_\chigh$ is tangent to
the path constant at $\chigh$, and $\tilo_d$ is tangent to
the path constant at $d$; therefore 
$\llangle\tilo^{(D\duer A)}_\chigh,\ \tilo_d \rrangle_{TA}$ is
tangent to the path constant at $\langle\chigh, d\rangle_A$, and 
is therefore zero. For the second term, $\Bar{\phi}^{(D\duer A)}$ 
is tangent to the path $t\timesCst\phi$ and $\Bar{d}_1^B$ is tangent
to the path $t\timesB d_1$, so 
$\llangle \Bar{\phi}^{(D\duer A)},\Bar{d}_1^B \rrangle_{TA}$ is
tangent to the path $\langle t\timesCst\phi, t\timesB d_1\rangle_A$
and by (v) of \ref{prop:dudvb} this is $t\langle\phi, d_1\rangle_A.$
\rule{0pt}{12pt}Altogether we have that 
$\llangle\mathcal{X}, \xi\rrangle_{TA} = \langle\phi, d_1\rangle_A.$
Using (\ref{eq:rz1}), we have that the RHS of (\ref{eq:rz4}) is
$$
\langle\phi, d_1\rangle_A - \langle \psi, d_1\rangle_B
$$
and this is equal to the RHS of (\ref{eq:rz3}) 
by the proof of \ref{thm:dualduality}. 

The rest of the proof is now straightforward. 
\ked

For a single vector bundle $E\to M$, the pairing of $E^*$ with its 
dual $E^{**}$ can be identified in a straightforward way with the 
pairing of $E$ with its dual. For double vector bundles it is first 
necessary to ensure that pairings are chosen in a consistent way. 

Consider a double vector bundle $(D;A,B;M)$ and assign signs to
the two upper structures as in Figure~\ref{fig:corner}(a) in order to
show that we pair the duals according to (\ref{eq:3duals}). 
Now, referring to Figure~\ref{fig:proof2}(b), 
\begin{figure}[htb]
\setlength{\unitlength}{1cm}
\begin{picture}(5,4)
\put(0,2){$\begin{matrix}
	&&\ominus&&\\
        &D&\mlra&B&\\
        &&&&\\
\oplus  &\Big\downarrow&&\Big\downarrow&\\
        &&&&\\
        &A&\mlra&M&\\
        &&&&\\
        &&\mbox{(a)}&&\\
\end{matrix}$}
\put(4.5,2){$\begin{matrix}
                      && \oplus                 &\\
                      &D\duer A&\longrightarrow &C^*\\
                      &&&\\
               \ominus&\Big\downarrow& &\Big\downarrow\\
                      &&&\\
                      &A&\longrightarrow &M\\
                      &&&\\
 		      &&\mbox{(b)}&\\
\end{matrix}$
}
\put(9,2){$\begin{matrix}
                      &&  \oplus                &\\
                      &D\duer B&\longrightarrow &B\\
                      &&&\\
           \ominus    &\Big\downarrow& &\Big\downarrow\\
                      &&&\\
                      &C^*&\longrightarrow &M\\
                      &&&\\
 		      &&\mbox{(c)}&\\
\end{matrix}$
}
\end{picture}
\caption{\ \label{fig:corner}} 
\end{figure}
we assign signs in such a way that each of $A, B, C^*$ has one 
positive and one negative arrow approaching it, and each of
$D$, $D\duer A$, $D\duer B$ has one positive and one negative
arrow departing from it; see 
Figure~\ref{fig:corner}(b)(c). We therefore, for example, 
take the pairing of the duals of $D\duer A$ to be 
\begin{equation}
\label{eq:DduerA}
\lhangle D, D\duer A\duer C^*\rhangle_B
= \langle D\duer A, D\duer A\duer C^*\rangle_{C^*}
	- \langle D, D\duer A\rangle_A.
\end{equation}
\begin{prop}
For the isomorphism $Z_B\co D\duer B\to (D\duer A\duer C^*)$, 
$$
\lhangle\Phi, \Psi\rhangle_{C^*} 
	= \langle\Phi, Z_B(\Psi)\rangle_{C^*},\qquad
\langle d, \Psi\rangle_{B} = -\lhangle d, \Z_B(\Psi)\rhangle_B
$$
for compatible $\Phi\in D\duer A$, $\Psi\in D\duer B$, $d\in D$. 
\end{prop}

\pf
The first is the definition of $Z_B$. The second follows from
$Z_B\duer C^* = Z_A$.
\ked

If we insert these equations into (\ref{eq:DduerA}), we get 
(\ref{eq:3duals}). Thus the signing on $D\duer A$ is compatible
with that on $D$. 

For ordinary vector bundles $E_1$ and $E_2$ on the same base $M$, 
one could take the view that a pairing of $E_1$ with $E_2$ is what
enables one to construct a cotangent double vector bundle with
sides $E_1$ and $E_2$. In a similar way, three double vector bundles
with suitably overlapping sides can be completed to a cotangent
triple vector bundle if and only if any two of them are the duals
of the third. 

\begin{df}
Consider three double vector bundles as in 
Figure~{\rm \ref{fig:cornering}(a)}, together with six pairings:
$\langle~,~\rangle_{E_1}$ of $D_2$ and $D_3$ over $E_1$, 
$\langle~,~\rangle_{E_2}$ of $D_3$ and $D_1$ over $E_2$, 
$\langle~,~\rangle_{E_3}$ of $D_1$ and $D_2$ over $E_3$, and
$\langle~,~\rangle_{1}$ of $C_1$ and $E_1$ over $M$, 
$\langle~,~\rangle_{2}$ of $C_2$ and $E_2$ over $M$, 
$\langle~,~\rangle_{3}$ of $C_3$ and $E_3$ over $M$, such that
each pairing of $D$ bundles is a pairing of double vector bundles
(as defined in {\rm \ref{prop:dudvb}}) with respect to the pairing 
of the relevant cores and sides. Then if
$$
\langle D_2, D_3\rangle_{E_1} = 
	\langle D_1, D_2\rangle_{E_3} - \langle D_1, D_3\rangle_{E_2}
$$
holds, we say that the system is a \emph{cornering} of $D_1$ with 
$D_2$ and $D_3$. 
\end{df}

\begin{figure}[htb]
\setlength{\unitlength}{1cm}
\begin{picture}(10,4.5)
\put(1,3){$\begin{matrix}
                      &&                        &\\
                      &&&D_3\\
                      &&&\\
                      && &\Big\downarrow\\
                      &&&\\
                      &D_1&\longrightarrow &E_2\\
            \end{matrix}
                      $}
\put(2.9, 3.5){\vector(2,-1){1}}                    
\put(1.6, 1.8){\vector(2,-1){1}}                    
\put(2.9, 1.8){\vector(2,-1){1}}                    

\put(2.7,1.7){$\begin{matrix}
                     &&                        &\\
                     &D_2&\longrightarrow &E_1\\
                     &&&\\
                     &\Big\downarrow&&\Big\downarrow\\
                     &&&\\
                     &E_3&\longrightarrow & M\\
		     &&&\\
		     &\mbox{(a)}&&\\	
               \end{matrix}
                     $}

\put(7,3){$\begin{matrix}
                      &&                        &\\
                      &T(T^*E)&\longrightarrow &T(E^*)\\
                      &&&\\
                      &\Big\downarrow& &\Big\downarrow\\
                      &&&\\
                      &TE&\longrightarrow &TM\\
            \end{matrix}
                      $}
\put(7.8, 3.5){\vector(2,-1){1.8}}                    
\put(9.4, 3.5){\vector(2,-1){1.8}}                    
\put(7.8, 1.8){\vector(2,-1){1.8}}                    
\put(9.4, 1.8){\vector(2,-1){1.8}}                    
\put(9.7,1.2){$\begin{matrix}
                     &&                        &\\
                     &T^*E&\longrightarrow &E^*\\
                     &&&\\
                     &\Big\downarrow&&\Big\downarrow\\
                     &&&\\
                     &E&\longrightarrow & M\\
		     &&&\\
		     &\mbox{(b)}&&\\	
               \end{matrix}
                     $}
\end{picture}
\caption{\ \label{fig:cornering}} 
\end{figure}

Clearly, choosing any double vector bundle in a cornering, the other 
two double bundles may be identified with its duals, and the 
cornering may be identified with the lower sides of the cotangent 
triple vector bundle associated with the chosen double. 

\begin{rmk}\rm
For an ordinary vector bundle $E$ one may form the cotangent
triple of $D = TE$. Now the canonical diffeomorphism
between $T^*TE$ and $TT^*E$ \citep{AbrahamM} is, since $E$ is
a vector bundle, an isomorphism of double vector bundles, and
so the triple $T^*TE$ is isomorphic to the tangent prolongation
of $T^*E$, as shown in Figure~\ref{fig:cornering}(b). 
Now the pairing of the bundles over $E$ in the left face gives 
rise to the canonical 1--form on $T^*E$, and the pairing of the 
bundles in the roof gives rise to the canonical 1--form on $T^*E^*$. 
\end{rmk}

\section{Duals of triple vector bundles}
\label{sect:dtvb}

Consider now a general triple vector bundle $\tvb$ as in 
Figure~\ref{fig:tvb}(a). Dualize $\tvb$ over the base $D_1$. Each of
the upper double vector bundles of which $\tvb\to D_1$ is a side
has a dual which is familiar from \S\ref{sect:proldual}. Following
the example of \ref{thm:T*D}, we complete the cube as in 
Figure~\ref{fig:dualsoftvbs}(a). 

\begin{figure}[htb]
\setlength{\unitlength}{1cm}
\begin{picture}(10,4)
\put(1,3){$\begin{matrix}
                      &&                        &\\
                      &\tvb\duer D_1&\longrightarrow &K_3\duer E_3\\
                      &&&\\
                      &\Big\downarrow& &\Big\downarrow\\
                      &&&\\
                      &D_1&\longrightarrow &E_3\\
            \end{matrix}
                      $}
\put(1.8, 3.5){\vector(2,-1){1.7}}                    
\put(3.3, 3.5){\vector(2,-1){1.7}}                    
\put(1.8, 1.8){\vector(2,-1){1.7}}                    
\put(3.3, 1.8){\vector(2,-1){1.7}}                    

\put(3.3,1.2){$\begin{matrix}
                     &&                        &\\
                     &K_2\duer E_2&\longrightarrow &W^*\\
                     &&&\\
                     &\Big\downarrow&&\Big\downarrow\\
                     &&&\\
                     &E_2&\longrightarrow & M\\
		     &&&\\
		     &\mbox{(a)}&&\\	
               \end{matrix}
                     $}

\put(7.2,3){$\begin{matrix}
                      &&                        &\\
                      &\tvb^{P_Y}&\longrightarrow &D_1\\
                      &&&\\
                      &\Big\downarrow& &\Big\downarrow\\
                      &&&\\
                      &D_2&\longrightarrow &E_3\\
            \end{matrix}
                      $}
\put(7.8, 3.5){\vector(2,-1){1.7}}                    
\put(9.1, 3.5){\vector(2,-1){1.7}}                    
\put(7.8, 1.8){\vector(2,-1){1.7}}                    
\put(9.1, 1.8){\vector(2,-1){1.7}}                    
\put(9.3,1.2){$\begin{matrix}
                     &&                        &\\
                     &D_3^P&\longrightarrow &E_2\\
                     &&&\\
                     &\Big\downarrow&&\Big\downarrow\\
                     &&&\\
                     &E_1&\longrightarrow & M\\
		     &&&\\
		     &\mbox{(b)}&&\\	
               \end{matrix}
                     $}
\end{picture}
\caption{\ \label{fig:dualsoftvbs}} 
\end{figure}

\begin{thm}
There is a triple vector bundle as shown in 
Figure~{\rm \ref{fig:dualsoftvbs}(a)} in which the four 
vertical sides are dual double vector bundles as just described. 
\end{thm}

We omit the details of this. As with the case $\tvb = TD$, five
of the six faces are double vector bundles of known types and
the main work is to show that the roof --- which belongs to a new
class of examples --- is a double vector bundle, and calculate
its core, which is $K_1\duer C_1$. The proof follows exactly 
the same outline as in \ref{prop:roofisdvb}, though steps
involving derivatives must be replaced with forms of the
interchange laws. 

Notice that in Figure~\ref{fig:dualsoftvbs}(a), two of the three 
upper double vector bundles are standard duals of the
double vector bundles in the corresponding positions in
Figure~\ref{fig:tvb}(a). Two of the lower double vector bundles
are duals of core double vector bundles of $\tvb$. 

The core double vector bundles of $\tvb\duer D_1$ are given in
Figure~\ref{fig:cdvbsdual}. The ultracore is $E_1^*$, the dual
of the bundle which in the original was diagonally opposite
$\tvb$ in the plane perpendicular to the axis of dualization.  

\begin{figure}[hb]
\setlength{\unitlength}{1cm}
\begin{picture}(10,3)
\put(0.3,1.4){$\begin{matrix}
                      &&                        &\\
                      &K_1\duer C_1&\longrightarrow &W^*\\
                      &&&\\
                      &\Big\downarrow& &\Big\downarrow\\
                      &&&\\
                      &C_1&\longrightarrow &M\\
                      &&&\\
		      &&\mbox{(a)}&\\
            \end{matrix}
                      $}
\put(4.8,1.41){$\begin{matrix}
                      &&                        &\\
                      &D_2\duer E_3&\longrightarrow &C_2^*\\
                      &&&\\
                      &\Big\downarrow& &\Big\downarrow\\
                      &&&\\
                      &E_3&\longrightarrow &M\\
                      &&&\\
		      &&\mbox{(b)}&\\
            \end{matrix}
                      $}
\put(9.3,1.4){$\begin{matrix}
                      &&                        &\\
                      &D_3\duer E_2&\longrightarrow &C_3^*\\
                      &&&\\
                      &\Big\downarrow& &\Big\downarrow\\
                      &&&\\
                      &E_2&\longrightarrow &M\\
                      &&&\\
		      &&\mbox{(c)}&\\
            \end{matrix}
                      $}
\end{picture}
\caption{\ \label{fig:cdvbsdual}} 
\end{figure}

The relationship between the three duals of $\tvb$ is 
embodied in the cotangent quaternary vector bundle of $\tvb$, 
as shown in Figure~\ref{fig:quat}. 

\begin{figure}[htb]
\setlength{\unitlength}{8mm} 
\begin{picture}(10,14)(-1.3,3)
\put(4,15){$\begin{matrix}
                      &&   &\\
                      &T^*\tvb&\longrightarrow &\tvb\duer D_2\\
                      &&&\\
                      &\Big\downarrow& &\Big\downarrow\\
                      &&&\\
                      &\tvb\duer D_1&\longrightarrow &K_3\duer E_3\\
            \end{matrix}
                      $}
\put(4.5,15.5){\vector(-1,-1){4}}        
\put(5.2,15.5){\vector(1,-1){4}}       
\put(6.5,15.5){\vector(-1,-1){4}}        
\put(7.3,15.5){\vector(1,-1){4}}       
\put(4.6,13.2){\vector(-1,-1){4}}        
\put(5.2,13.2){\vector(1,-1){4}}       
\put(6.5,13.2){\vector(-1,-1){4}}        
\put(7.3,13.2){\vector(1,-1){4}}       
\put(0.1,10.2){$\begin{matrix}
                      &&   &\\
                      &\tvb&\longrightarrow &D_2\\
                      &&&\\
                      &\Big\downarrow& &\Big\downarrow\\
                      &&&\\
                      &D_1&\longrightarrow &E_3\\
            \end{matrix}
                      $}
\put(8.7,10.2){$\begin{matrix}
                      &&   &\\
                      &\tvb\duer D_3&\longrightarrow &K_1\duer E_1\\
                      &&&\\
                      &\Big\downarrow& &\Big\downarrow\\
                      &&&\\
                      &K_2\duer E_2&\longrightarrow &W^*\\
            \end{matrix}
                      $}
\put(0.9,10.9){\vector(1,-1){4.2}}        
\put(2.2,10.9){\vector(1,-1){4.2}}       
\put(0.8,8.2){\vector(1,-1){4}}        
\put(2.2,8.2){\vector(1,-1){4}}       
\put(9.7,10.9){\vector(-1,-1){4.2}}        
\put(11.5,10.9){\vector(-1,-1){4.2}}       
\put(9.5,8.4){\vector(-1,-1){4}}        
\put(11.5,8.4){\vector(-1,-1){4}}       
\put(4.8,5.2){$\begin{matrix}
                      &&   &\\
                      &D_3&\longrightarrow &E_1\\
                      &&&\\
                      &\Big\downarrow& &\Big\downarrow\\
                      &&&\\
                      &E_2&\longrightarrow &M\\
            \end{matrix}
                      $}
\end{picture}
\caption{\ \label{fig:quat}}
\end{figure}

Denote dualization of $\tvb$ along the three axes by $X$, $Y$ and 
$Z$. In terms of Figure~\ref{fig:tvb}(a), take $Z$ to be dualization 
along the vertical axis, $Y$ to be along $D_3$ 
and $X$ to be along $D_2$. 
Compositions such as $ZXZ$ are triple versions of the operation
$P$ studied in \S\ref{sect:ddvb}. Precisely, 
applying $ZXZ$ to $\tvb$ in Figure~\ref{fig:tvb}(a) applies $P$ to
the rear face and to the front face; denote this by $P_Y$. This
operation may also be regarded as reflection of $\tvb$ in the plane
through $\tvb$, $D_3$, $M$ and $E_3$; see 
Figure~\ref{fig:dualsoftvbs}(b). Notice that each face has
been flipped in the sense that it cannot be returned to its original
position by a proper rotation of the cube. Further, the core double 
vector bundle which lies in the plane through $\tvb$, $D_3$, $M$ and 
$E_3$ is left fixed; the other two are flipped and interchanged. 

With similar definitions of $P_X$ and $P_Z$ we have
\begin{equation}
\label{eq:PPP}
P_X = YZY = ZYZ,\quad
P_Y = ZXZ = XZX, \quad
P_Z = XYX = YXY, 
\end{equation}
each of $P_X$, $P_Y$, $P_Z$ having order 2. Equivalently, each of
$XY$, $YZ$, $ZX$ has order~3. 

New in the triple case are the products $Q_Z = ZXYZ$, 
$Q_X = XYZX$ and $Q_Y = YZXY$ and their inverses. 
It is easily found from (\ref{eq:PPP}) that $Q_X$, $Q_Y$ and
$Q_Z$ have order~3 and that $Q_ZQ_YQ_X = I$. Curiously, the
equation $Q_XQ_YQ_Z = I$ or, equivalently, $(XYZ)^4 = I$, is not 
a consequence of (\ref{eq:PPP}), but it may be verified directly 
by calculating the effect on $\tvb$. We now have:
 
\begin{thm}
\label{thm:vb3}
The group of operations on $\tvb$ generated by 
$X$, $Y$ and $Z$ satisfies the relations 
$X^2 = I$, $Y^2 = I$, $Z^2 = I$, $(XYZ)^4 = I$, $(YZX)^4 = I$, 
$(ZXY)^4 = I$, together with {\rm (\ref{eq:PPP})}. 
\end{thm}

By a calculation with \citep{GAP4}, the group defined by these 
relations has order~72. Denote the group of operations 
generated by $X$, $Y$ and $Z$ by $\mathcal{VB}_3$. It is 
straightforward to find more than 36 distinct elements of 
$\mathcal{VB}_3$ and so it must have order 72. It thus cannot 
be, as one might have expected, a subgroup of the full symmetry 
group of the hypercube, which has order~384 \citep{Coxeter}. 

This shows that the situation with double vector bundles, in which
the operations generated by dualization can be identified with
symmetries of the cotangent triple, does not extend in the 
analogous fashion to triple vector bundles and symmetries of
the hypercube. 

\section{General principles}
\label{sect:gp}

On the basis of the duality theory for duals and triples, we
may formulate some likely principles for the duality of general
multiple vector bundles. It may be that the proofs are mainly 
a matter of acquiring sufficient motivation and notation, but
we cannot rule out the possibility that new phenomena arise 
with increasing dimension.  
 
There are three groups associated with an $n$--fold vector bundle
$\mathcal{N}$. Firstly, the various operations generated by flips
of the constituent double vector bundles generate an action of 
the symmetric group $\mathcal{S}_n$. Secondly, there is an obvious 
sense in which individual $n$--fold vector bundles may have more 
symmetry than others, as we remarked in the case of $T^*D$ 
and $TD$ in \S\ref{sect:tvbs}.

Of most interest, however, is the group $\mathcal{VB}_n$ generated
by the dualization operations. We have seen that $\mathcal{VB}_2$
is $\mathcal{S}_3$ and that $\mathcal{VB}_3$ has order 72. Further, 
the subgroup of $\mathcal{VB}_3$ generated by $XYXZ$, $YZYX$ and 
$ZXZY$ has order 12 and is normal, with quotient isomorphic to 
$\mathcal{S}_3$. 

In the general case, the $n$ duals of an $n$--fold vector bundle 
$\mathcal{N}$ and $\mathcal{N}$ itself form the $(n+1)$ lower 
$n$--faces of an $(n+1)$--fold vector bundle, which may be 
completed to be the cotangent $(n+1)$--fold vector bundle of 
$\mathcal{N}$, or of any of the duals of $\mathcal{N}$. 

The $(n+1)$ upper $n$--faces of $T^*\mathcal{N}$ consist of the
$n$ cotangent $n$--fold vector bundles of the upper $(n-1)$--faces
of $\mathcal{N}$, together with one $n$--fold vector bundle of a 
new type, which incorporates data from all of the structure 
of $\mathcal{N}$. It is reasonable to conjecture that if
$X_1, \dots, X_n$ denote the dualization operations, each of 
order 2, then we have, for each $1\leq k\leq n$ and each 
string $i_1, i_2, \dots, i_k$ of $k$ distinct elements of 
$\{1, \dots, n\}$, 
$$
(X_{i_1}X_{i_2}\cdots X_{i_k})^{k+1} = 1.
$$

  \def\cprime{$'$}

\end{document}